\documentclass[12pt]{article}%
\usepackage{amsfonts}
\usepackage{fancyhdr}
\usepackage{tabularx}
\usepackage{mathtools}
\usepackage[colorlinks = True]{hyperref}
\hypersetup{
    colorlinks = true,
    linkcolor = blue,
    anchorcolor = blue,
    citecolor = blue,
    filecolor = blue,
    urlcolor = blue
    }
\usepackage[a4paper, top=2.5cm, bottom=2.5cm, left=1.5cm, right=1.5cm]%
{geometry}
\usepackage{times}
\usepackage[table,xcdraw]{xcolor}
\usepackage{amsmath}
\usepackage{amsthm}
\usepackage[capitalize, nameinlink]{cleveref}
\usepackage{changepage}
\usepackage{enumitem}
\usepackage{mathrsfs}
\usepackage{tcolorbox}
\usepackage{stackrel,amssymb}
\usepackage{graphicx}

\usepackage{appendix}
\usepackage[backend=biber,
style=numeric,
maxnames=50
]{biblatex}
\def\hmath$#1${\texorpdfstring{{\rmfamily\textit{#1}}}{#1}}

\newtheorem{theorem}{Theorem}
\newtheorem{mainresult}[theorem]{Main Result}
\newtheorem{corollary}[theorem]{Corollary}
\newtheorem{observation}[theorem]{Observation}
\newtheorem{conjecture}[theorem]{Conjecture}
\newtheorem{definition}[theorem]{Definition}
\newtheorem{lemma}[theorem]{Lemma}
\newtheorem{remark}[theorem]{Remark}

\newtheorem{proposition}[theorem]{Proposition}
\newtheorem{problem}[theorem]{Problem}

\numberwithin{theorem}{subsection}

\newcommand{\anagram}{{\mathsf{A}}}
\newcommand{\afl}{{\mathcal{AFL}}}
\newcommand{\poset}{{\mathcal{P}}}
\newcommand{\canonicalword}{{\mathsf{cw}}}
\newcommand{\derangement}{\mathsf{D}}
\newcommand{\anagraphs}{\mathcal{AG}}
\newcommand{\derangementgraph}{\mathcal{DG}}
\def\hmath$#1${\texorpdfstring{{\rmfamily\textit{#1}}}{#1}}

\author{Kiril Bangachev}
\title{Enumerative and Structural Aspects Of\\
Anagrams Without Fixed Letters}
\date{}
\begin{document}
\maketitle
\abstract{
 For the word $\omega = \underbrace{11\ldots 1}_{x_1}\underbrace{22\ldots2}_{x_2}\ldots\underbrace{nn\ldots n}_{x_n},$ denote 
by $\mathsf{A}(x_1, x_2, \ldots, x_n)$ the number of its 
anagrams without fixed letters. While the function $\mathsf{A}()$ bears significant importance to economic theory \cite{MCKELVEY1997411}, it is not known whether it can be computed in polynomial time. The desire to answer efficiently certain queries related to this function motivates our study of its combinatorial properties. Our first main result shows that $\mathsf{A}(x_1, x_2, \ldots, x_n)\pmod{p}$  can be efficiently computed for any prime $p = O((\log n)^{1/3}).$ Our second main result establishes that the function $\mathsf{A}()$ is Schur-concave, which means that certain ordinal queries about $\mathsf{A}()$ can be answered in linearithmic time.
\\

Our second direction of study is structural. We introduce the anagraph, which generalizes derangement graphs. For $(x_1, x_2, \ldots, x_n)\in \mathbb{Z}_{\ge 0}^n,$ $\mathcal{AG}(x_1, x_2, \ldots, x_n)$ is a graph on vertex set all words over the alphabet $[n]$ which have exactly $x_i$ letters $i.$ Two vertices are adjacent if they are anagrams without fixed letters of each other. Our main result fully determines the $n$-tuples $(x_1, \ldots, x_n)$ for which the anagraph is connected and leads to a linear algorithm for this task.  We end with a conjecture, which fits into the ongoing debate about the connection between hamiltonicity and vertex-transitivity \cite{transitive}.
\\

One contribution of the current paper is a systematic development of techniques for analyzing anagrams without fixed letters. We illustrate the power of these techniques with further arithmetic, ordinal, and structural results.
}\\

\vspace*{1cm}

\vspace*{1cm}
\noindent
 \textsc{MIT Electrical Engineering $\&$
Computer Science Department, 
    Cambridge, MA, 02139}\par\nopagebreak
    \noindent
  \textit{Email address:} \href{mailto:kirilb@mit.edu}{kirilb@mit.edu}

\newpage
\section{Introduction}
\subsection{Background}
In the classical question of counting derangements, one considers permutations without fixed elements. That is, permutations $\sigma$ of $[n] = \{1,2,\ldots, n\}$ which satisfy that $\sigma(i)\neq i$ for all $i \in [n].$ This enumerative question is very well-understood. If $d(n)$ is the number of such permutations, the following well-known identities, among others, hold.

\begin{theorem}
\label{thm:derangement}
For the number of derangements $d(n),$ the following hold.
\begin{enumerate}
    \item $d(n)=(n-1)(d(n-1)+d(n-2))$ when $n \ge 2.$
    \item $\displaystyle d(n)=\left[ \frac{n!}{e}+\frac{1}{2}\right].$
    \item $\displaystyle d(n)=n!\sum_{i=0}^n\frac{(-1)^i}{i!}.$
\end{enumerate}
\end{theorem}

\noindent
While enumerating derangements is a simple and well-understood question, the problem motivates many different directions of active research today. These directions are both enumerative and structural.

\subsubsection*{Enumerative Directions}
Most important to the current paper is the  generalization of enumerating anagrams without fixed letters. It appears in a wide range of contexts such as combinatorics \cite{laguerre}, theoretical computer science \cite{tmne}, and economic theory \cite{MCKELVEY1997411}. 
The setup is the following.  For a word $\omega$ over the alphabet $[n],$ one is interested in anagrams \footnote{An anagram of a word $\omega$ is another word $\omega'$ which has the exact same letters, counted with multiplicities, but in a potentially different order. For example, the anagrams of $\omega = 121$ are $112,121,$ and $211.$}  $\omega'$ of $\omega$ such that $\omega'$ and $\omega$ have different letters at each position. We denote the number of such anagrams $\omega'$ by $\anagram(\omega).$
For example, if $\omega = 1123,$ then its anagrams without fixed letters are only $2311$ and $3211,$ so the desired number is $\anagram(\omega) =2.$
Clearly,  this number does not depend on the order of letters in $\omega,$ but only depends on the number of their individual appearances. For that reason, if $\omega$ has exactly $x_i$ letters $i,$ now on we will simply write $\anagram(x_1, x_2, \ldots, x_n)$ instead of $\anagram(\omega).$ Over the years, different properties of the function $\anagram : \mathbb{Z}_{\ge 0}^n \longrightarrow\mathbb{Z}_{\ge 0}$ have been studied such as its asymptotics when $x_1 = x_2 = \ldots = x_n$ \cite{asymptotic} and its relationship to bounding the number of totally mixed Nash equilibria \cite{tmne}. We discuss these in more depth in the next section.\\

\noindent
A very closely related 
generalization to anagrams
without fixed letters is the following. We phrase it in the language of Christmas presents used by Penrice \cite{asymptotic}. 
Suppose that there are $n$ disjoint groups of people $X_1, X_2, \ldots, X_n,$ where $|X_i| = x_i.$ 
The members of these groups want to exchange Christmas presents such that every person gives and receives exactly one present and no person receives a present from a member of their own group (including themselves). We denote the number of ways to accomplish this by
$\derangement(x_1, x_2, \ldots, x_n).$ When $x_1 = x_2 = \ldots = x_n =1,$ clearly $\derangement(x_1, x_2, \ldots, x_n) = d(n).$ One can also note that $\derangement(x_1, x_2, \ldots, x_n)= x_1!x_2!\cdots x_n!\anagram(x_1, x_2, \ldots, x_n),$ which means that the two functions have almost the same properties. Nevertheless, as will become apparent especially in \cref{sec:order}, sometimes it is easier to analyze the function $\derangement()$ to argue about $\anagram()$ and, for that reason, we also introduce the function $\derangement().$ For a quick demonstration of the convenience of $\derangement(),$ note that 

\begin{equation}
\label{eq:permanent}
\derangement(x_1, x_2, = \ldots = x_n) = 
per\begin{pmatrix}
M_1 \mathbf{1}\;\ldots \mathbf{1}\; &\\
\mathbf{1}\; M_2 \ldots \mathbf{1}\; &\\
\; \; \;\ddots\\
\mathbf{1}\; \mathbf{1}\;  \ldots M_n &\\
\end{pmatrix}.
\end{equation}
On the diagonal there are $n$ zero matrices $M_i,$ where $M_i \in \mathbb{Z}^{x_i\times x_i}.$ The rest of the entries are equal to 1.\\

\noindent
Many other generalizations of derangements also appear in literature. Different examples of these can be found, for instance, in \cite{fderangements}, \cite{rderangements}, and \cite{subgroupderangements}.

\subsubsection*{Structural Directions}
In addition to enumerative, structural properties of derangements and their generalizations are also studied. Specifically, of interest is the graph $\derangementgraph(n)$ which has as its vertices all permutations $\mathcal{S}_n$ and two nodes are connected if and only if one is a permutation without fixed elements of the other. That is, $\sigma, \tau \in V(\derangementgraph(n))$ are connected whenever $\sigma(i)\neq \tau(i)$ for all $i \in [n].$ Different properties of this graph have been studied such as hamiltonicity, path-hamiltonicity, edge-pancyclicity, its spectrum, and its independence number. For an overview of these, see \cite{edgepancyclic} and the references cited in it.\\

\noindent
In the current paper, we introduce and analyze the following natural generalization of the graph $\derangementgraph(n)$ in the context of anagrams without fixed letters.

\begin{definition} For an $n$-tuple $(x_1, x_2, \ldots, x_n)\in \mathbb{Z}^n_{\ge 0},$ we define the anagraph $\anagraphs(x_1, x_2, \ldots, x_n)$ as the following graph. Its nodes are all the words over the alphabet $[n]$ which have exactly $x_i$ letters $i$ for all $i \in [n].$ Two words $\omega$ and $\omega'$ are connected if and only if they differ at every position.
\end{definition}
\subsection{Previous Work}
\subsubsection*{Enumerative Questions}
While classical derangements are easy to enumerate as shown in \cref{thm:derangement} (and, clearly, extremely efficient algorithms for finding $d(n)$ exist), it is an open question whether one can find the number $\anagram(x_1, x_2, \ldots, x_n)$ with a polynomial-time algorithm even for the equinumerous case \linebreak $x_1 = x_2 = \cdots = x_n = k$ \cite{tmne}. The permanent formula \cref{eq:permanent} does not help as finding permanents of general $0\slash1$ matrices is an $\#$P-hard problem \cite{VALIANT1979189}. The following theorem, however, demonstrates that understanding $\anagram(x_1, x_2, \ldots, x_n)$ is an important, even if difficult, task. 

\begin{theorem}[\cite{MCKELVEY1997411,tmne}]
\label{thm:gametheory}
Consider a game with $n$ players, where player $i$ has $m_i\ge 1$ options to choose their action from. Then, a sharp upper bound on the number of totally mixed Nash equilibria in the game is
$\anagram(m_1 - 1, m_2 - 1, \ldots, m_n - 1).$
\footnote{A \textit{totally mixed} Nash equilibrium is a Nash equilibrium in which every player chooses every action available to them with positive probability.}
\end{theorem}

\noindent
Even though it is an open question how to simply compute the number $\anagram(\underbrace{k,k,\ldots,k}_n),$ the asymptotic growth of $\anagram()$ in this equinumerous case has been established in \cite{asymptotic}\footnote{In the original paper, the statement is $\lim_{n\longrightarrow +\infty}
\frac{\derangement(\overbrace{k,k,\ldots,k}^n)}{(nk)!} = e^{-k},$ but this is clearly equivalent to the statement in \cref{thm:asymptotic}. We prefer to phrase the result in terms of the function $\anagram()$ as anagrams without fixed letters are the main object of study in the current paper.}.

\begin{theorem}[\cite{asymptotic}]
\label{thm:asymptotic}
For any fixed $k\in \mathbb{N},$ one has 
$$
\lim_{n\longrightarrow +\infty}
{\anagram(\underbrace{k,k,\ldots,k}_n)}{\binom{{nk}}{\underbrace{k,k,\ldots,k}_n}^{-1}} = e^{-k}.
$$
\end{theorem}
\noindent
This theorem generalizes the second statement in \cref{thm:derangement}. Informally, it says that asymptotically an $e^{-k}$ fraction of all words containing exactly $k$ times each letter in $[n]$ are anagrams without a fixed letter of the word  $\omega = \underbrace{11\ldots1}_k\underbrace{22\ldots2}_k\ldots\underbrace{nn\ldots n}_k.$\\

\noindent
Another line of research exploits the ingenious connection between the function $\derangement()$ and Laguerre polynomials discovered by Gillis and Evans
\cite{laguerre}. This connection can be used, for example, to elegantly derive expressions for $\derangement()$ when the alphabet size (i.e., the number $n$) is small \cite{laguerre}.\\

\noindent
The importance of understanding  $\anagram(),$ justified by \cref{thm:gametheory}, motivates the study of exact (non-asymptotic) properties of this function. Our combinatorial analysis of such properties leads to efficient (sometimes even nearly-linear!) algorithms that answer certain queries about the function $\anagram().$

\subsubsection*{Structural Questions}
The graph $\derangementgraph(n)$ turns out to have a surprisingly rich structure. Particularly important to the current paper are the following two theorems. Before stating them, however, we recall a few common graph properties. 

\begin{definition}
\label{def:hamiltonicity}
A simple graph $G$ on $m$ vertices is
\begin{itemize}
    \item \textit{Hamiltionian} if there exists a simple cycle of length $m.$
    \item \textit{Hamilton-connected} if every pair of distinct vertices is joined by a simple path of length $m.$
    \item \textit{Pancyclic} if there exists a simple cycle of length $k$ for any $3\le k \le m.$
    \item \textit{Edge-Pancyclic} if there exists a simple cycle of length $k$ containing the edge $e$ of $G$ for any $e$ and  $3\le k \le m.$
\end{itemize}
\end{definition}

\noindent
It turns out that $\derangementgraph(n)$ satisfies all of these extremely strong properties. 

\begin{theorem}[\cite{edgepancyclic}]
\label{thm:edgepancyclic}
The graph $\derangementgraph(n)$ is edge-pancyclic for all $n \ge 4.$
\end{theorem}

\begin{theorem}[\cite{hamiltonconnected}]
\label{thm:hamiltonconnected}
The graph $\derangementgraph(n)$ is hamilton-connected for all $n \ge 4.$
\end{theorem}

\noindent
Naturally, we are curious to establish generalizations of these results for the graph $\anagraphs(x_1, x_2, \ldots, x_n).$ The end goal is finding necessary and sufficient conditions for an $n$-tuple $(x_1, x_2, \ldots, x_n)$ which guarantee that the corresponding anagraph is hamiltonian/path-hamiltonian/pancyclic/edge-pancyclic. It turns out, however, that some non-trivial work is needed even to determine when an anagraph is connected as will become apparent in \cref{se:anagraph}.\\

\noindent
The study of hamiltonicity of anagraphs is also of interest due to the ongoing efforts to understand 
which
connected vertex-transitive graphs are hamiltonian (see \cite{transitive} and references in it for a discussion of this connection). It can be easily shown that any anagraph is vertex-transitive (see \cref{obs:vtransitiveanagraph}) and we determine the condition for connectivity of anagraphs in \cref{thm:generalconnectivity}.\\

\noindent
The question of determining hamiltonicity (and related properties) of an anagraph  is also of interest from a computer-science point of view. It is well known that it is NP-hard to verify all of the properties listed in \cref{def:hamiltonicity} for a general graph $G$ \cite{reducability}. Even more, brute-force algorithms are particularly inefficient for anagraphs as anagraphs can have exponentially many vertices and edges in terms of their description as follows from \cref{obs:sizeobservation}.

\subsection{Main Questions and Results}
\subsubsection*{Enumerative Aspects of Anagrams Without Fixed Letters}
Since the number $\anagram(x_1, x_2, \ldots, x_n)$ has been elusive from a computational point of view so far, we focus on designing efficient algorithms that determine different properties of it. Our results show that the function $\anagram()$
has a surprisingly simple and elegant arithmetic and ordinal behaviour.\\

\noindent
Our first direction of study deals with the number-theoretic properties of $\anagram()$.

\begin{problem}
\label{prob:arithmetic}
Given an $n$-tuple $(x_1, x_2, \ldots, x_n)$ and a positive integer $m,$ determine the residue of $\anagram(x_1, x_2, \ldots, x_n)$ when divided by $m.$ Can this be done with an efficient algorithm?
\end{problem}

\noindent
We provide the following partial answer to this question.

\begin{mainresult}
\label{mnrslt:prime} There exists a polynomial-time algorithm that computes 
$A(x_1, \ldots, x_n)\pmod{p}$ for any $n$-tuple $(x_1, \ldots, x_n)\in \{0,1,\ldots, M\}^n$ and prime number $p = O((\log n + \log \log M)^{1/3}).$  In the special case $x_1 = \cdots = x_n = k,$ the result can be improved to primes of order $O((\log n + \log \log k)^{1/2}).$
\end{mainresult}

\noindent
Our proof of this result reduces computing  $A(x_1, x_2, \ldots, x_n)\pmod{p}$ to a few computations of the same type, but of much smaller size. The small size allows us to use exponential-time algorithms on them. Our result becomes especially elegant in the case $p = 2.$ In it, one simply needs to xor all the binary representations of $x_i$ and then check whether the resulting vector is non-zero. Namely, we have the following theorem (it is stated more precisely in \cref{cor:parity}). 

\begin{theorem} 
\label{thm:paritynew}
The number $\anagram(x_1,x_2,\ldots,x_n)$ is odd if and only if in the binary representations of $x_1,x_2,\ldots,x_n$ on every position there is an even number of ones.
\end{theorem}

\noindent
We present our arithmetic results in \cref{sec:arithmetic}. There, we first prove two other number-theoretic statements about $\anagram(),$
\cref{thm:primemodulo,thm:tail}. Not only are these results useful when proving \cref{mnrslt:prime}, but they also illustrate more simply the techniques used in \cref{mnrslt:prime}.\\

\noindent
In addressing the ordinal properties of $\anagram(),$ we pose the following question.

\begin{problem}
\label{prob:ordinal}
Given two $n$-tuples $(x_1, x_2, \ldots, x_n)$ and $(y_1, y_2, \ldots, y_n)$ of non-negative integers, determine which of the numbers $\anagram(x_1, x_2, \ldots, x_n)$ and 
$\anagram(y_1, y_2, \ldots, y_n)$ is larger.
\end{problem}

\noindent
We study this question in two different regimes depending on whether the sums $\displaystyle \sum_i x_i$ and $\displaystyle\sum_i y_i$ are equal. In the equality case, our main result relies on a specific poset defined over the $n$-tuples of non-negative integers with a fixed sum. The ordering relation is given as follows.

\begin{definition}
\label{def:majorize}
For two $n$-tuples of non-negative integers $(x_1, x_2, \ldots, x_n)$ and $(y_1, y_2, \ldots, y_n)$ such that $\displaystyle \sum_{i  =1}^n x_i = \sum_{i  =1}^n y_i,$ we say that $(y_1, y_2, \ldots, y_n)$ majorizes $(x_1, x_2, \ldots, x_n)$ and write\linebreak $(x_1, x_2, \ldots, x_n)\preceq (y_1, y_2, \ldots, y_n)$ if the following condition holds. For two permutations $\pi$ and $\sigma$ such that 
$x_{\pi(1)}\ge x_{\pi(2)}\ge \cdots \ge x_{\pi(n)}$ and 
$y_{\sigma(1)}\ge y_{\sigma(2)}\ge \cdots \ge y_{\sigma(n)},$ it is the case that 
$\displaystyle \sum_{j = 1}^i y_{\sigma(j)}\ge 
\sum_{j = 1}^i x_{\pi(j)}
$ holds for all $i \in [n].$
\end{definition}

\noindent
This ordering relation appears in many analytic and combinatorial inequalities, most well-known of which are the inequality of Karamata \cite{karamata} and the Schur-convexity property \cite{schurconcave}. In the current paper, it shows up in the following theorem.

\begin{mainresult}[Schur-Concavity of $\anagram()$ and $\derangement()$]
\label{thm:karamatatype}
Suppose that $(x_1, x_2, \ldots, x_n)\preceq (y_1, y_2, \ldots, y_n)$ in the sense of \cref{def:majorize}. Then, 
\begin{equation*}
    \begin{split}
       & \anagram(x_1, x_2, \ldots, x_n)\ge 
\anagram(y_1, y_2, \ldots, y_n), \; \; \; \text{ and }\\
 & \derangement(x_1, x_2, \ldots, x_n)\ge 
\derangement(y_1, y_2, \ldots, y_n).
    \end{split}
\end{equation*}
\end{mainresult}

\noindent
In other words, $\anagram()$ and $\derangement()$ are Schur-concave in the sense of \cite{schurconcave}.

\begin{corollary} There exists an $O(n\log n)$ algorithm, which determines which of the numbers\linebreak $\anagram(x_1, x_2, \ldots, x_n)$ and $\anagram(y_1, y_2, \ldots, y_n)$ is greater, provided that the $n$-tuples $(x_1, x_2, \ldots, x_n)$ and $(y_1, y_2, \ldots, y_n)$ are comparable with respect to $\succeq.$
\end{corollary}

\noindent
The same technique used in the proof of \cref{thm:karamatatype} also allows us to derive \cref{thm:varsumder,thm:varsuman}, which handle certain special cases in the regime $\displaystyle \sum_{i  =1}^n x_i \neq \sum_{i  =1}^n y_i.$ More specifically, in this regime we analyze the minimal possible difference, i.e. $\displaystyle \sum_{i  =1}^n x_i =1+\sum_{i  =1}^n y_i.$

\subsubsection*{Structural Aspects of Anagrams Without Fixed Letters}

\noindent
So far, we have discussed only enumerative questions about anagrams without fixed letters. These correspond to an analysis of an extremely local property - the degree - of an anagraph as will become apparent in \cref{obs:sizeobservation}. However, derangement graphs have a very rich global structure as demonstrated in \cref{thm:hamiltonconnected,thm:edgepancyclic}. This motivates us to study global properties of anagraphs. Our first question relates to the, perhaps, most prominent global property - connectivity.

\begin{problem} Determine the $n$-tuples $(x_1, x_2, \ldots, x_n)\in \mathbb{Z}_{\ge 0}^n$ for which the graph $\anagraphs(x_1, x_2, \ldots, x_n)$ is connected.
\end{problem}

\noindent
We fully resolve this question with the following theorem.

\begin{mainresult}
\label{thm:firstconnected}
The anagraph $\anagraphs(x_1, x_2, \ldots, x_n)$ is connected if and only if one of the following conditions is satisfied:
\begin{enumerate}
    \item $n = 1.$
    \item $n = 2$ and $x_1 = x_2 = 1.$
    \item $n \ge 3, $ the inequality $\max(x_1, x_2, \ldots, x_n) <\frac{1}{2}\sum_{j = 1}^n x_j$ holds, and $(x_1, x_2, \ldots, x_n)\neq (1,1,1).$
\end{enumerate}
\end{mainresult}

\noindent
Leaving some edge-cases aside, connectivity of the anagraph is equivalent to 
$\max(x_1, x_2, \ldots, x_n) <\frac{1}{2}\sum_{j = 1}^n x_j.$ On the other hand, one can simply observe (see \cref{obs:positivity} and \cref{obs:sizeobservation}) that the anagraph has an empty edgeset if and only if $\max(x_1, x_2, \ldots, x_n) >\frac{1}{2}\sum_{j = 1}^n x_j.$ This shows a very sharp transition between absence of edges and connectivity.\\

\noindent
Again, \cref{thm:firstconnected} leads to an efficient algorithm. 

\begin{corollary} There is an algorithm which determines whether $\anagraphs(x_1, x_2, \ldots, x_n)$ is connected on input $(x_1, x_2, \ldots, x_n)$ in time $O(n).$
\end{corollary}

\noindent
Of course, more interesting is the question of determining the properties in \cref{def:hamiltonicity}.

\begin{problem}
\label{prob:hamiltonicityproblem}
Determine the $n$-tuples $(x_1, x_2, \ldots, x_n)$ for which the graph $\anagraphs(x_1, x_2, \ldots, x_n)$ is hamiltonian (respectively, hamiltonian-connected, pancyclic, and edge-pancyclic).
\end{problem}

\noindent
We conjecture that every connected anagraph (with potentially finitely many exceptions) is at least hamiltonian in light of the curious connection between vertex-transitivity and hamiltonicity \cite{transitive}. In fact, only four connected vertex-transitive graphs on at least 3 vertices that are not hamiltonian are known \cite{transitive}!\\

\noindent
While we do not fully resolve \cref{prob:hamiltonicityproblem}, we show in 
\cref{prop:reduction}
that a reduction of the alphabet size - which is our main tool in proving \cref{thm:firstconnected} - also applies to it. This means that if our conjecture is true, it is enough to prove that every connected anagraph on an alphabet of size 4 is at least hamiltonian.

\section{Preliminaries}
\subsection{Further Notation and Terminology}
\subsubsection*{Words and Anagrams}
Throughout, we will denote by 
$\canonicalword(x_1, x_2, \ldots, x_n)$
the ``canonical word'' on alphabet $[n]$ with exactly $x_i$ letters $i,$ which is given by 
$\underbrace{11\ldots 1}_{x_1}\underbrace{22\ldots 2}_{x_2}\ldots\underbrace{nn\ldots n}_{x_n}.$ Most of the arguments in \cref{sec:arithmetic,sec:order} will be formulated for this word. If $\chi = \canonicalword(x_1, x_2, \ldots, x_n),$ for any other word $\omega$ over $[n]$ which has exactly $x_i$ letters $i,$ we will often decompose $\omega$ as 
$$
\omega = \chi_1(\omega)\chi_2(\omega)\cdots\chi_n(\omega),
$$
where each $\chi_i(\omega)$ is a subword of $\omega$ that contains exactly $x_i$ consecutive letters.\\

\noindent
For a word $\omega$ with $m$ letters, we will enumerate the positions by the positive integers from left to right. For example, when $\omega = 57911$ (so, $m = 5$) the word has letter 5 at position $1, $ letter $7$ at position 2, letter $9$ at position $3,$ and letter $1$ at positions $4$ and 5. For any $S\subseteq[m],$ we will denote by $\omega|_S,$ the word $\omega$ restricted to positions $S.$ For example, if $\omega = 57911,$ and $S = \{2,3,4\},$ then $\omega|_S = 791.$\\

\noindent
For a word $\omega,$ we will denote by $\afl(\omega)$ its set of anagrams without fixed letters. Clearly,\linebreak $|\afl(\omega)| = \anagram(\omega).$

\subsubsection*{The Majorization Relation}
We end this section with a few notes on the majorization relation defined in \cref{def:majorize}. First, we note that the definition can easily be extended to $n$-tuples of different length by simply adding zeros. Adding zeros to any $n$-tuple $(x_1, x_2, \ldots, x_n)$ can be done freely anywhere in the paper.

\begin{definition}
\label{def:majorizationposet}
For a positive integer $m,$ we define the following poset $\poset_m.$ Its elements are all positive integer sequences with sum $m$ and the poset relation given by $\preceq$ in \cref{def:majorize}.
\end{definition}

\noindent
We present a simple illustration of how the poset structure defined above appears in combinatorial inequalities. We will need the following inequality later on.

\begin{observation}
\label{obs:schurconvexityfactorials}
If $(y_1, y_2, \ldots, y_n)$ majorizes $(x_1, x_2, \ldots, x_n),$ then 
$y_1!y_2!\ldots y_n!\ge x_1!x_2!\ldots x_n!.$
\end{observation}
\begin{proof}
This fact follows simply from log-convexity of the gamma function \cite{gammaconvexity} and Karamata inequality \cite{karamata}. Nevertheless, we present a different proof technique, which will be useful when discussing \cref{thm:fixedsumanar}.\\

\noindent
First, note that whenever $t\ge s>0,$ it is the case that $(t+1)!(s-1)!\ge 
t!s!.
$
Now, assume without loss of generality that $y_1 \ge y_2 \ge \ldots\ge y_n$ and $x_1\ge x_2 \ge \ldots \ge x_n.$
Define the $n$-tuple\linebreak $(v^{(0)}_1, v^{(0)}_2, \ldots, v^{(0)}_n) =(y_1, y_2, \ldots, y_n).$ While $(v^{(k)}_1, v^{(k)}_2, \ldots, v^{(k)}_n)\neq (x_1, x_2, \ldots, x_n),$ do the following procedure. We will prove that it is well defined. 
\begin{enumerate}
    \item Find the smallest index $i$ such that $v^{(k)}_i \neq x_i.$
    \item Find the smallest index $j>i$ such that $v^{(k)}_{j-1} <v^{(k)}_i.$
    \item Update $v^{(k+1)}_{j-1} = v^{(k)}_{j-1} - 1, v^{(k+1)}_j = v^{(k)}_j + 1,$ $v^{(k+1)}_r = v^{(k)}_r$ for $r\not \in \{j,j-1\}.$
\end{enumerate}
First, note that after every update, $\sum_s v^{(k)}_s$ remains unchanged, so the $n$-tuples $(v^{(k)}_1, v^{(k)}_2, \ldots, v^{(k)}_n)$ are elements of the poset $\poset_{y_1 + y_2 + \cdots+y_n}.$ We show by induction that $(v^{(k)}_1, v^{(k)}_2, \ldots, v^{(k)}_n)\succeq (x_1, x_2, \ldots, x_n)$ holds for all $k.$ Indeed, this is true for $k = 0.$ Now, if it is true for some $k$ and $(v^{(k)}_1, v^{(k)}_2, \ldots, v^{(k)}_n)\neq (x_1, x_2, \ldots, x_n),$ we first need to show that steps 1 and 2 can be executed. Let $i$ be the minimal index such that $v^{(k)}_i \neq x_i.$ Clearly, as 
$(v^{(k)}_1, v^{(k)}_2, \ldots, v^{(k)}_n)\succeq (x_1, x_2, \ldots, x_n),$  it must be the case that\linebreak $v^{(k)}_i > x_i.$ Now, we claim that there exists some $j >i$ such that $v^{(k)}_j  < v^{(k)}_i.$ Indeed, if this were not true, one would derive the following contradiction.
$$
\sum_{u = 1}^n v^{(k)}_u  = 
\sum_{u = 1}^i v^{(k)}_i  + (n-i)v^{(k)}_i >
\sum_{u = 1}^i x_i + (n-i)x_i \ge 
\sum_{u = 1}^n x_u.
$$
Let $j$ be the minimal index such that $v^{(k)}_j <v^{(k)}_i.$ Clearly, $v^{(k)}_{j -1} = v^{(k)}_i >x_i \ge  x_{j-1}$ must hold. Thus, we can perform the update in step 3. One can very easily check that the new $n$-tuple\linebreak $(v^{(k+1)}_1, v^{(k+1)}_2, \ldots, v^{(k+1)}_n)$ also majorizes $(x_1, x_2, \ldots, x_n)$ by the choice of $i$ and $j.$\\
Therefore, the procedure indeed terminates with $(x_1, x_2, \ldots, x_n).$ The statement follows from the update rule 3. and the observation that whenever $t\ge s>0,$ it is the case that $(t+1)!(s-1)!\ge 
t!s!.
$
\end{proof}

\subsection{Simple Enumerative Observations}
We first note that the functions $\anagram $ and $\derangement$ are symmetric. That is, for any permutation $\pi,$\linebreak $\anagram(x_1, x_2, \ldots, x_n) = \anagram(x_{\pi(1)}, x_{\pi(2)}, \ldots, x_{\pi(n)})$ and similarly for $\derangement.$ We continue with determining for which words $\omega,$ there exists anagrams without fixed letters.

\begin{observation}[\cite{tmne}]
\label{obs:positivity}
$\anagram(x_1, x_2, \ldots, x_n)>0$ (and 
$\derangement(x_1, x_2, \ldots, x_n)>0$) is equivalent to \linebreak
$\displaystyle \max(x_1, x_2, \ldots, x_n)\le 
\frac{1}{2}\sum^n_{i = 1} x_i.
$
\end{observation}
\begin{proof}
Without loss of generality, let $x_1 \ge x_2\ge \cdots\ge x_n.$\\
First, suppose that $x_1 >
\frac{1}{2}\sum^n_{i = 1} x_i.
$ Then, clearly, any anagram of $\canonicalword(x_1, x_2, \ldots, x_n)$ will have at least one letter 1 among its leftmost $x_i$ positions, so there are no anagrams without fixed letters.\\
On the other hand, if $x_1 \le \frac{1}{2}\sum^n_{i = 1} x_i,$ one can easily check that the cyclic permutation with $x_1$ positions to the left $\omega'$ is an anagram without fixed letters of $\omega.$ More precisely, we take\linebreak
 $\omega' = \underbrace{22\ldots2}_{x_2}\underbrace{33\ldots3}_{x_3}\ldots,\underbrace{nn\ldots n}_{x_n}\underbrace{11\ldots1}_{x_1}.$
\end{proof}

\noindent
We continue with two very simple observations about the arithmetic and ordinal structures of $\anagram()$ as a prelude to our main results.
To illustrate a simple property of the arithmetics of $\anagram,$ which will also be useful later on, we make the following simple observation about derangements.

\begin{observation} $\anagram(\underbrace{1,1,\ldots,1}_n)\equiv n-1 \pmod{2}.$
\end{observation}
\begin{proof}
This follows easily as $\anagram(\underbrace{1,1,\ldots,1}_n) = d(n),$ $d(0) = 1, d(0) = 1,$ and the first relation in \cref{thm:derangement}.
\end{proof}

\noindent
We will vastly improve the above observation in \cref{cor:parity}. We also make a simple observation about the ordinal structure of $\anagram().$

\begin{observation}Two $n$-tuples $(x_1, x_2, \ldots, x_n)$ and 
$(b_1, b_2, \ldots, b_n)$ of non-negative integers are given such that $\max(b_1, b_2, \ldots, b_n)\le \frac{1}{2}\sum_{i = 1}^n b_i.$ Then,
$$
\anagram(x_1, x_2, \ldots, x_n)\le 
\anagram(x_1+b_1, x_2+b_2, \ldots, x_n+b_n).
$$
\end{observation}
\begin{proof} Consider the words $\omega_1 = \canonicalword(b_1, b_2, \ldots, b_n),$ $\omega_2 = \canonicalword(x_1, x_2, \ldots, x_n)$ and 
$\omega = \omega_1\omega_2.
$ We will simply construct an injection $f:\afl(\omega_2)\longrightarrow\afl(\omega).$ Since $\max(b_1, b_2, \ldots, b_n)\le \frac{1}{2}\sum_{i = 1}^n b_i,$ the word $\omega_1$ has an anagram without fixed letters $\omega_1'.$ Then, $f$ takes the simple form $f(\chi) = \omega'_1\chi$ for any $\chi \in \afl(\omega_2).$
\end{proof}

\noindent
We continue with computing two very simple, specific, cases of $\anagram(x_1, x_2, \ldots, x_n).$

\begin{observation}[\cite{tmne}]
\label{obs:halfletter1}
If $x_1  = x_2 + \cdots +x_n,$ then 
$\anagram(x_1, x_2, \ldots, x_n) = \binom{x_1}{x_2, \ldots, x_n}.$
\end{observation}
\begin{proof}
Note that the anagrams without fixed letters of $\canonicalword(x_1, x_2, \ldots, x_n),$ are exactly the words of the form $\chi \underbrace{1,1,\ldots, 1}_{x_1}$ where $\chi$ is a word over $[n]\backslash\{1\}$ having exactly $x_i$ letters $i$ for all $i>1.$
\end{proof}

\begin{observation}
\label{obs:evaltwooneone}
$\anagram(2,\underbrace{1, 1, \ldots, 1}_{n-2}) = 
\frac{1}{2}(d(n) - 2d(n-1)-d(n)).
$
\end{observation}
\begin{proof}
We will prove that $\derangement(2,\underbrace{1, 1, \ldots, 1}_{n-2}) = d(n) - 2d(n-1)-d(n),$ which is enough. Let $X_1 = \{1,2\}$ and $X_i = \{i\}$ for $i \ge 3.$ We want to find the number of permutations $\sigma: \bigcup_{i}X_i\longrightarrow \bigcup_{i}X_i$ such that $\sigma(X_i)\cap X_i = \emptyset$ for all $i.$ First, note that any such permutation is necessarily a derangement of $\omega = 12\ldots n.$ Now, we will count how many derangements of $\omega$ do not satisfy the condition\linebreak $\sigma(X_i)\cap X_i = \emptyset$ for all $i.$ This condition can be violated by three types of derangement's.
\begin{itemize}
    \item When $\sigma(1) = 2, \sigma(2) = 1.$ There are clearly $d(n-2)$ such derangements.
    \item When $\sigma(1) = 2, \sigma(2) \neq 1.$ There are clearly $d(n-1)$ such derangements.
    \item When $\sigma(1) \neq  2, \sigma(2) = 1.$ There are clearly $d(n-1)$ such derangements.
\end{itemize}
In total, we count $d(n) - 2d(n-1) - d(n-2)$ such derangements.
\end{proof}

\subsection{Simple Observations About the Anagraph}
First, we state two simple observation about the size of the anagraph without proof.

\begin{observation}
\label{obs:sizeobservation}
For a $n$-tuple $(x_1, x_2, \ldots, x_n),$ one has:
\begin{enumerate}
    \item $\displaystyle|V(\anagraphs(x_1, x_2, \ldots, x_n))| = \binom{x_1 + x_2 + \cdots + x_n}{x_1, x_2, \ldots, x_n}.$ 
    \item For every $v\in V(\anagraphs(x_1, x_2, \ldots, x_n)),$ it is the case that 
    $deg(v) = \anagram(x_1, x_2, \ldots, x_n).$
\end{enumerate}
\end{observation}

\noindent
Another, very important property of anagraphs, is their high symmetry. 

\begin{observation} 
\label{obs:vtransitiveanagraph}
For any $n$-tuple $(x_1, x_2, \ldots, x_n),$ the anagraph $\anagraphs(x_1, x_2, \ldots, x_n)$ is\linebreak vertex-transitive.
\end{observation}
\begin{proof}
Let $\lambda$ and $\omega$ be two words in $V(\anagraphs(x_1, x_2, \ldots, x_n)).$ Since each letter $i$ appears the same number of times in $\lambda$ and $\omega,$ there exists a permutation $\sigma$ (of positions of letters) such that $\sigma(\lambda) = \omega.$ Clearly, $\sigma,$ when extended to $V(\anagraphs(x_1, x_2, \ldots, x_n)),$ is a graph isomorphism.
\end{proof}
\subsection{Approaches to Anagrams Without Fixed Letters}
One of our main contributions is developing techniques for a systematic study of the enumerative and structural aspects of anagrams without fixed letters. The main approaches that we use in the paper are three:
\begin{itemize}
    \item \textbf{Equivalence Classes}. This technique is used when discussing the arithmetic of $\anagram()$ in \cref{thm:mainprimemodthm,thm:primemodulo,thm:tail}.
    We split $\afl(\chi)$ for a word $\chi$ into simple equivalence classes, the size of which we can easily analyze. This technique
    can also be phrased in terms of group actions. For example, in \cref{thm:mainprimemodthm,thm:primemodulo}, the proof is equivalent to studying the orbits of a group action of $\mathcal{S}_{x_1}\times \mathcal{S}_{x_2}\times \cdots \times \mathcal{S}_{x_n}$ on $\afl(\canonicalword(x_1, x_2, \ldots, x_n)).$ In \cref{thm:tail}, we utilize a group action of $\mathcal{C}_m$ on 
    $\afl(\canonicalword(\underbrace{k,\ldots, k}_m, x_{m+1}, \ldots, x_n)).$
    \item \textbf{Recurrence Relations.} Specifically, we analyze a recurrence relation for $\derangement$ and $\anagram,$ which reduces simultaneously the alphabet size and the word length (see \cref{thm:recurrence}). While such a recurrence relation does not yield an efficient algorithm for computing $\anagram(x_1, x_2, \ldots, x_n),$ it turns out to be extremely useful for proving inequalities like \cref{thm:fixedsumder,thm:varsumder,thm:varsuman}. The main insight is that it is relatively easy to control and compare coefficients  with which the values of $\derangement()$ for ``smaller'' words appear in \cref{thm:recurrence}.
    \item \textbf{Alphabet Reductions.} An alphabet reduction is part of our technique with recurrence relations. However, it appears more explicitly in \cref{prop:reduction}. We show that under a specific reduction of the alphabet size, properties such as hamiltonicity and connectivity are preserved. In the case of connectivity, this allows us to only consider the setting of $n \le 4$ in order to determine when an anagraph is connected (see \cref{thm:generalconnectivity}).
\end{itemize}
\section{Arithmetic Properties}
\label{sec:arithmetic}
We begin the section on arithmetic with two properties, which  are simpler to prove than \cref{mnrslt:prime}. Nevertheless, both the statements and techniques used in them are useful later on.
\subsection{The Case of Equinumerous Letters of Prime Order}
\label{sec:prime}
\begin{theorem}
\label{thm:primemodulo}
For any prime number $p$ and positive integer $n$, the following congruence holds.
$$\anagram(\underbrace{p,p,\ldots,p}_{n})\equiv \anagram(\underbrace{1,1,\ldots,1}_{n})\pmod{p^3}.$$
\end{theorem}
\noindent
Before we prove this theorem, we note that it is stronger than \cref{thm:paritynew} when $p=2.$ Indeed, \cref{thm:paritynew} only gives the congruence modulo $2,$ but \cref{thm:primemodulo} gives it modulo $8 = 2^3.$ The case of $p = 2$ also illustrates that the third power in \cref{thm:primemodulo} is optimal. When $p = 2, n = 3,$ we have that $\anagram(2,2,2) = 10$ and $\anagram(1,1,1) = 2$ (see \cite[p.141]{laguerre} for these values), so \linebreak
$\anagram(2,2,2)\not \equiv \anagram(1,1,1) \pmod{2^4}.$
\begin{proof} Let $\chi = \canonicalword(\underbrace{p,p,\ldots,p}_{n}).$ For any $\omega\in\afl(\chi,)$ where
$\omega = \chi_1(\omega)\chi_2(\omega)\cdots\chi_n(\omega),$ call the \textit{subword} $\chi_i(\omega)$ of $\omega$ a \textit{good subword} if it contains at least 2 distinct letters. Clearly, any word $\omega$ has either 0 or at least 2 good subwords. For each index $i$ and letter $j,$ denote by $\alpha_{i,j}(\omega)$ the number of letters $j$ in $\chi_i(\omega).$ Using the numbers $\alpha_{i,j},$ we introduce the following equivalence relation over $\afl(\chi)$. For $\omega, \lambda\in \afl(\chi)$ we have $\omega \sim \lambda$ if and only if $\alpha_{i,j}(\omega) = \alpha_{i,j}(\lambda)$ for all $i,j.$ 

\noindent
Note that the number of good subwords is a well-defined function over the defined equivalence classes. This allows us to analyze equivalence classes based on the number of good subwords they contain. We consider three cases.\\

\noindent
\textbf{Case 1.} Classes that contain $0$ good subwords. Note that each of these equivalence classes contains a single word of the form 
$$\underbrace{\sigma(1)\sigma(1)\ldots\sigma(1)}_p
  \underbrace{\sigma(2)\sigma(2)\ldots\sigma(2)}_p\ldots
  \underbrace{\sigma(n)\sigma(n)\ldots\sigma(n)}_p,  
$$
where $\sigma$ is a permutation of $[n]$ for which $\sigma(i)\neq i$ for all values of $i.$ Trivially, the number of such permutations is $d(n) = \anagram(\underbrace{1,1,\ldots,1}_{n}).$\\

\noindent
\textbf{Case 2.} Classes that contain exactly $2$ good subwords. We will explicitly count the number of words that appear in such classes and show that it is always divisible by $p^3.$\\
First, note that any word $\omega$ that has exactly two good subwords is of the following form. There exist two special letters $i$ and $j$ and indices $\pi(i),\pi(j)\not \in \{i,j\}$ such that $\chi_{\pi(i)}(\omega)$ and $\chi_{\pi(j)}(\omega)$ are both composed only of letters $i$ and $j.$ Any other $\chi_r(\omega)$ is composed of a single letter, different from $r,i,$ and $j.$ Suppose that for $k\in [n]\backslash\{i,j\},$ the subword containing (only) the letter $k$ is indexed by $\pi(k),$ i.e. $\chi_{\pi(k)} = \underbrace{kk\ldots k}_n.$
Say that in $\chi_{\pi(i)}(\omega)$ there are $1\le x\le p-1$ letters $i$ and $p-x$ letters $j.$ Therefore, in $\chi_{\pi(j)}(\omega)$ there are $p-x$ letters $i$ and $x$ letters $j.$\\

\noindent
We are ready to begin the count.
\begin{itemize}
    \item First, there are $\binom{n}{2}$ ways to choose the pair $(i,j).$
    \item Then, there are $\anagram(2,\underbrace{1,1,\ldots, 1}_{n-2})$ ways to choose $\pi(k)$ for $k \in [n].$ This is the case since we need to choose them in such a way that $\pi(i)\not \in \{i,j\}, \pi(j)\not \in \{i,j\},$ and $\pi(k)\neq k$ for $k \not \in \{i,j\}$ and any such choice satisfies the conditions.
    \item There are $\binom{p}{x}$ ways to arrange the letters $i$ and $j$ in $\chi_{\pi(i)}(\omega)$ and similarly for $\chi_{\pi(j)}(\omega).$
\end{itemize}

\noindent
Now, we simply sum over $x$ and conclude that the number of words which have exactly two good subwords is
$$
\sum_{x = 1}^{p-1}\binom{p}{x}^2\binom{n}{2}\anagram(2,\underbrace{1,1,\ldots, 1}_{n-2}).
$$
We consider three cases based on $p.$\\
\textbf{Case 2.1.} When $p>3.$ Then, 
$$
\sum_{x = 1}^{p-1}\binom{p}{x}^2 = \binom{2p}{p} - 2\equiv 0 \pmod{p^3},
$$
where the congruence modulo $p^3$ follows immediately from Wolstenholme's theorem \cite{wolstenholme}.\\

\noindent
\textbf{Case 2.2.} When $p = 3.$ Then, the above expression evaluates to $9n(n-1)\anagram(2,\underbrace{1,1,\ldots, 1}_{n-2}).$ When $n \equiv 0,1\pmod{3},$ clearly the expression is divisible by $27.$ When $n\equiv 2\pmod{3},$ we use the fact that $\anagram(2,\underbrace{1,1,\ldots, 1}_{n-2}) =
\frac{1}{2}(d(n) - 2d(n-1)-d(n-2))\equiv 0 \pmod{3}
$ (see \cref{obs:evaltwooneone}). The congruence follows simply from 
$d(n) = (n-1)(d(n-1) +d(n-2))$ and $d(0) = 1, d(1) = 0.$ \\

\noindent
\textbf{Case 2.3.} When $p =2.$ Then, the above expression evaluates to $2n(n-1)\anagram(2,\underbrace{1,1,\ldots, 1}_{n-2}).$ One can easily check from \cref{obs:evaltwooneone} that $\anagram(2,\underbrace{1,1,\ldots, 1}_{n-2})$ is even, from which the statement follows as $n(n-1)$ is also divisible by 2.\\

\noindent
\textbf{Case 3.} Classes that contain at least $3$ good subwords. Note that the number of words in the equivalence class of $\omega$ is exactly 
$$
\prod_{i = 1}^n\binom{p}{\alpha_{i,1}(\omega),
\alpha_{i,2}(\omega), \ldots, 
\alpha_{i,n}(\omega)
}. 
$$
Furthermore, whenever $\chi_i(\omega)$ is a good subword, the term
$\displaystyle \binom{p}{\alpha_{i,1}(\omega),
\alpha_{i,2}(\omega), \ldots, 
\alpha_{i,n}(\omega)
}$ is divisible by $p$ since all the terms $\alpha_{i,1}(\omega),
\alpha_{i,2}(\omega), \ldots, 
\alpha_{i,n}(\omega)$ are between $0$ and $p-1$ and $p$ is a prime. Therefore, for any word which has at least three good subwords, its equivalence class contains a number of words divisible by $p^3.$
\end{proof}

\noindent
\begin{remark}
\normalfont
While the statement of \cref{thm:primemodulo} might seem rather odd, it fits into a much broader family of results in elementary number theory of the form $f(p)\equiv f(1)\pmod{p^k}.$ When\linebreak $f(x) = n^x$ for some integer $n$ and $k = 1,$ this is Fermat's celebrated theorem. When $f(x) = \binom{2x}{x}$ and $k = 3,$ this is Wolstenholme's theorem \cite{wolstenholme}.
\end{remark}

\subsection{The Case of an Equinumerous Prefix} 
\label{sec:eqtail}
We continue our discussion of the arithmetic of the number of anagrams without fixed letters with the following theorem. Even though its proof is rather simple, the theorem has several interesting and perhaps surprising corollaries.

\begin{theorem}
\label{thm:tail}
Let $n\ge m >0,k>0$ and $x_{m+1}, x_{m+2}, \ldots, x_n$ be non-negative integers. Then
$$\anagram(\underbrace{k,k,\ldots,k}_{m},{x_{m+1},x_{m+2},\ldots, x_n})\equiv\anagram(\underbrace{k,k,\ldots,k}_{m}) \anagram({x_{m+1},x_{m+2},\ldots ,x_n})  \pmod{m}.$$
\end{theorem}
\noindent
We call this case the case of an equinumerous prefix as the $n$-tuple $\underbrace{k,k,\ldots,k}_{n},{x_{m+1},x_{m+2},\ldots, x_n}$ has the equinumerous prefix $\underbrace{k,k,\ldots,k}_{m}.$
\begin{proof}
Consider the word $\chi = \canonicalword(\underbrace{k,k,\ldots,k}_{m},{x_{m+1},x_{m+2},\ldots, x_n}).$ Now, for $i\in [m],$ define the functions $h_i$ as follows. For any $i$ and word $\lambda$ over $[n],$ the word $h_i(\lambda)$ is the same as $\lambda$ except that each letter $j\in [m]$ is replaced with $j + i$ if $j +i \le m$ and with $j+i -m$ if $j+i >m.$ Whenever $j \not \in [m],$ the letter $j$ remains unchanged. In other words, we cyclically permute the letters in the set $[m]$ with offset $i$ modulo $m$ and leave the letters in $[n]\backslash [m]$ unchanged.\\

\noindent
Now, one can easily check that if $\omega = \chi_1(\omega)\chi_2(\omega)\cdots\chi_{n}(\omega)$ is an anagram without fixed letters of $\chi,$ so 
is 
$$
f_i(\omega):=
h_i(\chi_{1 - i}(\omega))
h_i(\chi_{2 - i}(\omega))\cdots
h_i(\chi_{m - i}(\omega))
h_i(\chi_{m+1}(\omega))
h_i(\chi_{m+2}(\omega))
\cdots
h_i(\chi_{n}(\omega))
$$
for any $i \in [m],$
where $\chi_t(\omega):=\chi_{t+m}(\omega)$ whenever $t\le 0.$ Now, we define the following equivalence
relation over 
$\afl(\chi).$ For $\lambda, \omega\in \afl(\chi),$ it is the case that $\lambda\sim \omega$ if and only if there exists some $i\in [m]$ such that
$\lambda = f_i(\omega)$ (this is well defined, because if $\lambda = f_i(\omega),$ then 
$\omega = f_{m-i}(\lambda)$ when $i \neq m$ and $\lambda = \omega$ when $i = m$).\\

\noindent
Observe that if $\omega$ is such that its subword $\chi_{m+1}(\omega)\cdots \chi_n(\omega)$ contains at least one letter $i \in [m],$ the equivalence class of $\omega$ under $\sim$ contains exactly $m$ different words. As we are interested in the residue modulo $m,$ it is enough to count the anagrams without fixed letters of $\chi$ for which $\chi_{m+1}(\omega)\cdots \chi_n(\omega)$ only contains letters in $[n]\backslash [m].$ For any such word $\omega,$ it must also be the case that $\chi_{1}(\omega)\cdots \chi_m(\omega)$ only contains letters in $ [m].$ Therefore, the number of anagrams without fixed letters of $\chi$ that satisfy the desired property is exactly $\anagram(\underbrace{k,k,\ldots,k}_{m}) \anagram({x_{m+1},x_{m+2},\ldots ,x_n}),$
which completes the proof.
\end{proof}
\noindent
Now, we list without proof a few corollaries of \cref{thm:tail}.
\begin{corollary}
Let $k,m,$ and $n$ be positive integers. Then:
$$\anagram(\underbrace{k,k,\ldots,k}_{m},\underbrace{k,k,\ldots,k}_{n})\equiv \anagram(\underbrace{k,k,\ldots,k}_{m})\anagram(\underbrace{k,k,\ldots,k}_{n})\pmod{lcm(m,n)}.$$
\end{corollary}

\begin{corollary}
Let $n$ and $k$ be positive integers. If $m=nt+r$ for some non-negative $t$ and $r$, then $$\anagram(\underbrace{k,k,\ldots,k}_{m})\equiv \anagram(\underbrace{k,k,\ldots,k}_{r})( \anagram(\underbrace{k,k,\ldots,k}_{n}))^t \pmod{n}.$$
In particular, this means that whenever $m \equiv 1\pmod{n},$ one has $n | \anagram(\underbrace{k,k,\ldots,k}_{m}).$
\end{corollary}

\begin{corollary}
For two fixed positive integers $n$ and $k,$ define the sequence $a_m:=\anagram(\underbrace{k,k,\ldots,k}_{m}).$ Then $\displaystyle (a_i)_{i=1}^{+\infty}$ is eventually periodic modulo $n.$ 
\end{corollary}
\subsection{Main Result on Arithmetics}
We split the main result into the following theorems and propositions. \cref{mnrslt:prime} is a direct consequence of 
\cref{cor:algorithmprime}.

\begin{theorem}
\label{thm:mainprimemodthm}
Let $p$ be a prime and $(x_1, x_2, \ldots, x_n)$ be an $n$-tuple composed of non-negative integers less than $p^{N+1}.$ Let the base-$p$ representation of $x_i$ be $\overline{x_{i}^{(N)} x_{i}^{(N-1)}\ldots x_{i}^{(0)}}_{(p)}$ for each $i \in [n].$ Then, 
$$
\anagram(x_1, x_2, \ldots, x_n)\equiv
\prod_{j = 0}^N
\anagram(x_1^{(j)}, x^{(j)}_2, \ldots, x^{(j)}_n)\pmod{p}.
$$
\end{theorem}
\begin{proof}
For the word $\chi = \canonicalword(x_1, x_2, \ldots, x_n),$ define the functions $\alpha_{i,j}()$ over $\afl(\chi)$ and the equivalence relation $\sim$ in the same way as in the proof of \cref{thm:primemodulo}. 
Further, for any $\alpha_{i,j}(\omega)$, denote its base-$p$ representation by $\alpha_{i,j}(\omega) =
\overline{\alpha_{i,j}(\omega)^{(N)} \alpha_{i,j}(\omega)^{(N-1)}\ldots \alpha_{i,j}(\omega)^{(0)}}_{(p)}.
$
Now, observe that for any $\omega\in \afl(\chi),$ its equivalence class $[\omega]$ has size
$$
\prod_{i = 1}^n 
\binom{x_i}{\alpha_{i,1}(\omega),\alpha_{i,2}(\omega), \cdots,\alpha_{i,n}(\omega)}.
$$
Using Dickson's theorem about the residues of multinomial coefficients modulo prime numbers \cite{multinomial}, we conclude that the following congruence holds.
\begin{equation*}
    \begin{split}
        \prod_{i = 1}^n 
\binom{x_i}{\alpha_{i,1}(\omega),\alpha_{i,2}(\omega), \ldots,\alpha_{i,n}(\omega)} \equiv
 \prod_{i = 1}^n \prod_{s = 0}^{N}
\binom{x^{(s)}_i}{\alpha_{i,1}(\omega)^{(s)},\alpha_{i,2}(\omega)^{(s)}, \ldots,\alpha_{i,n}(\omega)^{(s)}} \pmod{p}.
    \end{split}
\end{equation*}
Now observe that for each $i\in [n],s\in \{0,1,\ldots,N\},$ the following two statements hold
\begin{equation}
    \label{eq:firstcongruence}
    \begin{split}
        \alpha_{i,1}(\omega)^{(s)}+\alpha_{i,2}(\omega)^{(s)}+ \cdots+\alpha_{i,n}(\omega)^{(s)}\equiv & 
        \; x^{(s)}_i \pmod{p},\\
        \alpha_{i,1}(\omega)^{(s)}+\alpha_{i,2}(\omega)^{(s)}+ \cdots+\alpha_{i,n}(\omega)^{(s)}\ge& 
        x^{(s)}_i,
    \end{split}
\end{equation}
and the multinomial coefficient $\binom{x^{(s)}_i}{\alpha_{i,1}(\omega)^{(s)},\alpha_{i,2}(\omega)^{(s)}, \ldots,\alpha_{i,n}(\omega)^{(s)}}$ is non-zero if and only if\linebreak $\alpha_{i,1}(\omega)^{(s)}+ \cdots+\alpha_{i,n}(\omega)^{(s)} = 
        x^{(s)}_i.$ Having this in mind, define the set
        $E(\chi)$ of \textit{no-carry-on} equivalence classes under $\sim$ as the set of equivalence classes for which 
        $\alpha_{i,1}(\omega)^{(s)}+ \cdots+\alpha_{i,n}(\omega)^{(s)} = 
        x^{(s)}_i$ holds for all $i$ and $s.$ It follows that 
$$
\anagram(x_1, x_2, \ldots, x_n)\equiv 
\sum_{[\omega]\in E(\chi)}
\prod_{i = 1}^n \prod_{s = 0}^{N}
\binom{x^{(s)}_i}{\alpha_{i,1}(\omega)^{(s)},\alpha_{i,2}(\omega)^{(s)}, \ldots,\alpha_{i,n}(\omega)^{(s)}} \pmod{p}.
$$
Our next task will be to better characterize $E(\chi).$ In particular, if $\chi^{(s)}:= 
\canonicalword(x_1^{(s)}, x_2^{(s)}, \ldots, x_n^{(s)})
$ for all $s\in \{0,1,\ldots, N\},$ we will show that 
$E(\chi)\cong E(\chi^{(0)})\times E(\chi^{(1)})\times \cdots\times E(\chi^{(N)}).$\\

\noindent
In order to proceed, we need some further notation. We partition the set of positions within each $\chi_i$ into sets with sizes equal to powers of $p.$ More precisely, for all triplets $1\le i \le n, 0 \le s \le N, 1\le k \le x_i^{(s)},$ we have that 
$$
T(\chi, i,s,k):=
\{
\sum_{u = 1}^{i-1}
x_u + 
\sum_{v = 0}^{s-1}p^vx_i^{(v)} + 
(k-1)p^s + 
r \; : \;
1 \le r \le p^s\}.
$$
Note that whenever $x_i^{(s)}$ is zero, the sets are not defined. We illustrate this definition to make it less abstract. For $p = 3,n = 2, x_1 = 2, x_2  = 7,$ and $\rho= \canonicalword(2,7),$ we have the following sets of positions:
$$
\canonicalword(2,7) = \rho = 
\underbracket{1}_{T(\rho, 1,0,1)}
\underbracket{1}_{T(\rho, 1,0,2)}
\underbracket{2}_{T(\rho, 2,0,1)}
\underbracket{222}_{T(\rho, 2,1,1)}
\underbracket{222}_{T(\rho, 2,1,2)}.
$$

\noindent
Now, observe that a class $e$ of $\afl(\chi)/\sim$ is no-carry-on if and only if there exists some $\lambda\in e$ such that for all $1\le i \le n, 0 \le s \le N, 1\le k \le x_i^{(s)},$ the subword $\lambda|_{T(\chi, i,s,k)}$ contains only a single letter (repeated $p^s$ times). Indeed, this follows directly from the definition of no-carry-on classes. 
For a no-carry-on class $e\in E(\chi),$ define $\lambda(e)$ to be alphabetically first word $\lambda\in e$ satisfying this property. Denote by $\ell_{e,i,s,k}$ the letter at positions $\lambda(e)|_{T(\chi, i,s,k)}.$\\

\noindent
Now, we construct injections in both ways between 
$E(\chi)$ and $E(\chi^{(0)})\times E(\chi^{(1)})\times \cdots\times E(\chi^{(N)}).$\\

\noindent
First, we construct an injective mapping 
$\mathsf{g}:E(\chi)\longrightarrow E(\chi^{(0)})\times E(\chi^{(1)})\times \cdots\times E(\chi^{(N)}).$
Take some $e\in E(\chi).$
For each $0\le s \le N,$ denote by 
$
\xi^{(s)}(e)
$ the following word:
$$\xi^{(s)}(e) = 
\ell_{e,1,s,1}\ell_{e,1,s,2}\cdots \ell_{e,1,s,x_1^{(s)}}
\ell_{e,2,s,1}\ell_{e,2,s,2}\cdots \ell_{e,2,s,x_2^{(s)}} \cdots
\ell_{e,n,s,1}\ell_{e,n,s,2}\cdots \ell_{e,n,s,x_n^{(s)}}.
$$
Then, 
$\mathsf{g}(e):= ([\xi^{(0)}], [\xi^{(1)}], \ldots,[\xi^{(N)}])$ is a well-defined injection.\\

\noindent
Second, we construct an injective mapping
$\mathsf{h}: E(\chi^{(0)})\times E(\chi^{(1)})\times \cdots\times E(\chi^{(N)})\longrightarrow E(\chi).$ Take an $(N+1)$-tuple $(e^{(0)}, e^{(1)}, \ldots, e^{(N)})\in
E(\chi^{(0)})\times E(\chi^{(1)})\times \cdots\times E(\chi^{(N)}).$ Then, $\mathsf{h}(e^{(0)}, e^{(1)}, \ldots, e^{(N)})$ is the word $\lambda\in \afl(\chi),$ such that $\lambda|_{T(\chi, i, s, k)}$ is composed of just $p^s$ times the letter $\lambda(e^{(s)})|_{T(\chi^{(s)},
i, 0, k)}.$ Again, clearly $\mathsf{h}$ is well-defined and injective. Furthermore, we can check that $\mathsf{g}\circ \mathsf{h}$ and 
$\mathsf{h}\circ \mathsf{g}$ both equal the identity.\\

\noindent
We make one further observation about $\mathsf{h}()$. Note that for each $i,j,s$ and $e^{(s)}\in E(\chi^{(s)}),$ the function 
$$(e^{(0)}, e^{(1)}, \ldots, e^{(N)})\longrightarrow
\alpha_{i,j}(\mathsf{h}
(e^{(0)}, e^{(1)}, \ldots, e^{(N)}))^{(s)}$$
depends only on $e^{(s)}$ and is invariant under changing the other $N$ coordinates $(e^{(t)})_{t\;:\; t\neq s}.$ Therefore, by abuse of notation, we will write 
$\alpha_{i,j}(e^{(s)})^{(s)}$ instead of $\alpha_{i,j}(\mathsf{h}
(e^{(0)}, e^{(1)}, \ldots, e^{(N)}))^{(s)}.$

\noindent
With this in mind, we can go back to computing $\anagram(x_1, x_2, \ldots, x_n)\pmod{p}.$ Namely, we have 
\begin{equation*}
    \begin{split}
&\anagram(x_1, x_2, \ldots, x_n)\equiv_p\\ 
&\sum_{[\omega]\in E(\chi)}
\prod_{i = 1}^n \prod_{s = 0}^{N} 
\binom{x^{(s)}_i}{\alpha_{i,1}(\omega)^{(s)},\alpha_{i,2}(\omega)^{(s)}, \ldots,\alpha_{i,n}(\omega)^{(s)}} = \\
&\sum_{e^{(0)} \in E(\chi^{(0)}),
e^{(1)} \in E(\chi^{(1)}),\ldots,
e^{(N)} \in E(\chi^{(N)})
}
\prod_{s = 0}^{N}
\prod_{i = 1}^n\\
&
\qquad \qquad
\binom{x^{(s)}_i}{\alpha_{i,1}(\mathsf{h}
(e^{(0)},  \ldots, e^{(N)}))^{(s)},
\alpha_{i,2}(\mathsf{h}
(e^{(0)},  \ldots, e^{(N)}))^{(s)}, \ldots,
\alpha_{i,n}(\mathsf{h}
(e^{(0)},  \ldots, e^{(N)}))^{(s)}} =\\
\end{split}
\end{equation*}

\begin{equation*}
    \begin{split}
&\sum_{e^{(0)} \in E(\chi^{(0)}),
e^{(1)} \in E(\chi^{(1)}),\ldots,
e^{(N)} \in E(\chi^{(N)})
}
\prod_{s = 0}^{N}
\prod_{i = 1}^n
\binom{x^{(s)}_i}{
\alpha_{i,1}(e^{(s)})^{(s)},
\alpha_{i,2}(e^{(s)})^{(s)}, \ldots,
\alpha_{i,n}(e^{(s)})^{(s)}} =\\
&\prod_{s = 0}^{N}
\sum_{e^{(s)}\in E(\chi^{(s)})}
\prod_{i = 1}^n
\binom{x^{(s)}_i}{
\alpha_{i,1}(e^{(s)})^{(s)},
\alpha_{i,2}(e^{(s)})^{(s)}, \ldots,
\alpha_{i,n}(e^{(s)})^{(s)}}.
    \end{split}
\end{equation*}
Now, using \cref{eq:firstcongruence} for $\chi^{(s)},$ we conclude that 
$$
\sum_{e^{(s)}\in E(\chi^{(s)})}
\prod_{i = 1}^n
\binom{x^{(s)}_i}{
\alpha_{i,1}(e^{(s)})^{(s)},
\alpha_{i,2}(e^{(s)})^{(s)}, \ldots,
\alpha_{i,n}(e^{(s)})^{(s)}}\equiv 
\anagram(
x_1^{(s)},
x_2^{(s)},
\cdots,
x_n^{(s)}
)\pmod{p},
$$
from which the statement follows.\end{proof}

\begin{proposition}
\label{prop:prefixreductionsprimemod}
Let $(t_1, t_2, \ldots, t_n)\in \{0,\ldots, p-1\}^n.$ Suppose that for each $c \in \{1,\ldots, p-1\},$ there are exactly $u_c$ numbers in $(t_1, t_2, \ldots, t_n)$ equal to $c.$ Suppose that $r_c\in \{0,\ldots, p-1\}$ is the residue of $u_c$ modulo $p.$ Then,
$$
\anagram(t_1, t_2, \ldots, t_n)\equiv 
\anagram(
\underbrace{1,\ldots, 1}_{r_1},
\underbrace{2,\ldots, 2}_{r_2},\ldots,
\underbrace{p-1,\ldots, p-1}_{r_{p-1}})
\prod_{c = 1}^{p-1}\anagram(\underbrace{c,c,\ldots, c}_{p})^{u_c-r_c}\pmod{p}.
$$
\end{proposition}
\begin{proof} We just repeatedly apply \cref{thm:tail} for $m = p.$
\end{proof}

\noindent
Now, we are ready to present our algorithmic result. 

\begin{corollary}
\label{cor:algorithmprime}
For a prime $p$ and an $n$-tuple 
$(x_1, x_2, \ldots, x_n)\in \{0,1,2,\ldots, M\}^{n},$ there exists an algorithm running in time 
$O(\exp(p^3) + n \times poly(\log M, p) + {p^2\log M})$ which determines $\anagram(x_1, x_2, \ldots, x_n)\pmod{p}.$ Furthermore, in the equinumerous case $x_1 = x_2 = \cdots = x_n = k,$ the running time reduces to \linebreak
$O(\exp(p^2) + n \times poly(\log k, p) + {p^2\log k}).$
\end{corollary}
\begin{proof}
Note that the input consists of $p, n, x_1, x_2, \ldots, x_n.$ Thus, its size is $O(\log p + n \log M).$ Our algorithm runs as follows:
\begin{enumerate}
    \item Compute the base-$p$ representations for each $x_i.$
    \item Compute the numbers $u_{c}^{(j)}$ and
    $r_{c}^{(j)}$ defined in \cref{prop:prefixreductionsprimemod} for each $j \in \{0,1,\ldots, N\}$ and $n$-tuple 
    $(x_1^{(j)}, x^{(j)}_2, \ldots, x^{(j)}_n).$ Then, compute $k_c$ defined as follows. 
    If $\sum_{j = 0}^N (u^{(j)}_c - r^{(j)}_c) = 0,$ then $k_c = 0.$ Otherwise, $k_c$ is defined as the number in $\{1,2,\ldots, p-1\}$ 
    having the same residue modulo $p-1$ as  $\sum_{j = 0}^N (u^{(j)}_c - r^{(j)}_c).$
    \item For each $j \in\{0,1,2,\ldots, N\}$ compute (if not already computed for a smaller $j$) and store the value of $\anagram(
\underbrace{1,\ldots, 1}_{r^{(j)}_1},
\underbrace{2,\ldots, 2}_{r^{(j)}_2},\ldots,
\underbrace{p-1,\ldots, p-1}_{r^{(j)}_{p-1}})\pmod{p}$ in a hash table. 
\item Compute  $\anagram(\underbrace{c,c,\ldots, c}_{p})\pmod{p}$ for all $c\in [p-1].$ 
    \item Compute and return the following residue modulo $p:$
    $$
    \prod_{j  = 0}^N
\anagram(
\underbrace{1,\ldots, 1}_{r^{(j)}_1},
\underbrace{2,\ldots, 2}_{r^{(j)}_2},\ldots,
\underbrace{p-1,\ldots, p-1}_{r^{(j)}_{p-1}})
\times
\prod_{c = 1}^{p-1}\anagram(\underbrace{c,c,\ldots, c}_{p})^{k_c}
    $$
\end{enumerate}
\noindent
Using \cref{thm:mainprimemodthm} and \cref{prop:prefixreductionsprimemod} together with the well known theorem due to Fermat stating that $x^p\equiv x \pmod{p}$ holds for all primes $p$ and integers $x,$ we conclude that the algorithm is correct. Now, we only need to argue about its running time.\\

\noindent
Step 1 can be clearly done in time $n \times poly(\log M, \log p).$ Now, note that each $x_i$ has $N = O(\log_{p}\log M)$ $= O(\frac{\log M}{\log p})$ digits in base $p.$ Therefore, step 2 can be performed in time\linebreak $Npn = O(n \times poly(p,\log M)).$ Step 3 can be performed in time $O(\exp (p^3))$ as follows. Using \cref{eq:permanent}, the complexity of finding $\anagram(
\underbrace{1,\ldots, 1}_{r^{(j)}_1},
\underbrace{2,\ldots, 2}_{r^{(j)}_2},\ldots,
\underbrace{p-1,\ldots, p-1}_{r^{(j)}_{p-1}})\pmod{p}$ is asymptotically the same
as finding the residue modulo ${p}$ of a 0/1 permanent of size $\sum_{i = 1}^{p-1}r_ii \le \frac{p^3}{2}.$ As proven in \cite{permanentmoduloprime}, this can be done in time $O(\exp{\frac{p^3}{2}}).$ Now, note that there are at most $p^{p-1}  = o(\exp(\frac{p^3}{2}))$
choices for $(r_1, r_2, \ldots, r_{p-1})\in \{0,1,\ldots, p-1\}^{p-1}.$ Thus, the number of new computations in Step 3 is $o(\exp(\frac{p^3}{2})),$ from which its total running time is $O(\exp({p^3})).$ In the equinumerous case, note that for each $j \in \{0,1,\ldots, N\}$ at most one of the numbers $r^{(j)}_c$ is non-zero and, thus, the size of the permanent is bounded by $\frac{p^2}{2}$ rather than $\frac{p^3}{2}.$
 Step 4 has the same analysis as step 3. In step 5, we simply need to perform $O(n + p\log p)$ multiplications of residues modulo $p,$ which can be done in time $n\times poly(p).$ 
\end{proof}

\begin{remark}
\normalfont
\label{rmk:primesize}
In particular,  for $p$ constant, $A(x_1, x_2, \ldots, x_n)\pmod{p}$ can be computed in linearithmic time. This shows that despite \cref{eq:permanent}, the problem of computing the function $\anagram()$ is - perhaps - significantly easier than computing the permanent of a general 0/1 matrix. It is well known that deciding whether the permanent of a general 0/1 matrix is divisible by 3 takes exponential time under the exponential time hypothesis \cite{permanentmodulo3}. More generally, \cref{cor:algorithmprime} shows that there exists a poly-time algorithm for computing $\anagram(x_1, x_2, \ldots, x_n)\pmod{p}$ when $p = O(({\log n + \log \log M})^{1/3})$ in the general case and when 
$p = O(({\log n + \log \log k})^{1/2})$ in the equinumerous case.
\end{remark}

\begin{remark}
\normalfont
\label{rmk:productofprimes}
One of our main motivations for studying $A(x_1, x_2, \ldots, x_n)\pmod{p}$ was that doing this computation for sufficiently many primes and, then using Chinese Remainder Theorem, will allows us to find $A(x_1, x_2, \ldots, x_n)\pmod{K}$ for some large $K.$ In the best case scenario, this $K$ would be so large that knowing $A(x_1, x_2, \ldots, x_n)\pmod{K}$ (and potentially an asymptotic result like \cref{thm:asymptotic}) will allow us to efficiently compute $A(x_1, x_2, \ldots, x_n)\pmod{K}.$ Nevertheless, a lot more work beyond primes of order $O(({\log n + \log \log M})^{1/3})$
is needed towards that direction. It is a well-known fact that 
$\prod_{q<x}q= e^{x(1 + o_x(1))}$ where the product is taken over primes less than $x$ \cite{chebyshevtheta}.
\end{remark}

\noindent
We end with restating the elegant characterization of the parity of $\anagram()$ in \cref{thm:paritynew} which follows directly from
\cref{thm:mainprimemodthm} and \cref{prop:prefixreductionsprimemod}. 

\begin{corollary}
\label{cor:parity}
 Let $x_1, x_2, \ldots, x_n$ be non-negative integers less than $2^{N+1}$ with binary representations $x_i = \overline{x_{i}^{(N)}, x_{i}^{(N-1)}, \ldots, x_{i}^{(0)}}_{(2)}.$ Then, 
the number $\anagram(x_1,x_2,\ldots,x_n)$ is odd if and only if the sum $\sum_{i = 1}^n{x_{i}^{(j)}}$ is even for all integers $j \in \{0,1,2\ldots, N\}.$ 
\end{corollary}

\section{Ordinal Properties}
\label{sec:order}
\subsection{Detour in a Recurrence Relation}
Perhaps surprisingly, the main tool in the proofs in this section is a recurrence relation. We choose to work with the function $\derangement()$ instead of $\anagram(),$ because, as we will see, the results we obtain for $\derangement()$ are actually stronger. To state the recurrence relation, we first need the following definition.

\begin{definition}
\label{def:recurencehelperfunction}
Define $f(x_1,x_2,\ell)$ for $x_1\ge x_2 \ge 0, \ell\ge 0$ as follows. $$f(x_1,x_2,\ell):=\sum_{\ell_1=0}^\ell\binom{x_1}{\ell_1}\binom{x_2}{\ell_1}(\ell_1)! \binom{x_1}{\ell-\ell_1}\binom{x_2}{\ell-\ell_1}(\ell-\ell_1)!.$$
\end{definition}

\noindent
Note that whenever $\ell > 2\min(x_1, x_2),$ it is the case that $f(x_1, x_2, \ell) = 0.$\\

\noindent
\cref{def:recurencehelperfunction} is motivated by the following proposition.
\begin{proposition}
\label{prop:finterpretation}
Let $X_1$ with $|X_1|=x_1$ and $X_2$ with $|X_2|=x_2$ be two disjoint sets. Then $f(x_1,x_2,\ell)$ is the number of ways to choose a subset $L\subseteq X_1\cup X_2$ of size  
$ \ell$  and construct an injective function $\sigma:L\rightarrow X_1\cup X_2$ which satisfies the following condition. No element of $X_1\cap L$ is mapped to an element in $X_1$ and, similarly, no element of $X_2\cap L$ is mapped to an element in $X_2.$
\end{proposition}
\begin{proof}
For fixed sizes $|X_1\cap L| =\ell_1 $ and $|X_2\cap L| = \ell-\ell_1,$ we can choose $X_1\cap L$ in $\binom{x_1}{\ell_1}$ ways and $\sigma(X_1\cap L)\subseteq X_2$ in  $\binom{x_2}{\ell_1}$ ways. For each choice, there are exactly $\ell_1!$ mappings from  $X_1\cap L$ to  $\sigma(X_1\cap L).$  We argue analogously for  $X_2\cap L$ and  $\sigma(X_2\cap L)$ and sum over $\ell_1.$ 
\end{proof}

\noindent
With the help of \cref{prop:finterpretation}, we can prove the following recurrence relation, which will be the main workhorse in the current section. 

\begin{theorem} For any $n$-tuple of non-negative integers $(x_1, x_2, \ldots x_n),$ the following equality holds.
\label{thm:recurrence}
$$\sum_{\ell=0}^{x_1+x_2}f(x_1,x_2,\ell) \derangement(x_1+x_2-\ell,x_3,\ldots,x_n)=\derangement(x_1,x_2,\ldots,x_n).$$
\end{theorem}
\begin{proof}
Let $X_1, X_2, \ldots, X_n$ be $n$ disjoint sets, where $|X_i| = x_i.$ We will count in two ways the number of bijections $\sigma$ of $X := \bigcup_{i = 1}^nX_i$ such that $X_i \cap \sigma(X_i) = \emptyset$ for all $i.$ By the definition of $\derangement(),$ this number is exactly $\derangement(x_1,x_2,\ldots,x_n).$\\

\noindent
For any $\sigma$ satisfying this property, denote $S(\sigma):=(X_1\cup X_2)\cap\sigma(X_1\cup X_2).$ Let $|S(\sigma)| = \ell(\sigma).$ Now, we will consider the following equivalence relation $\sim_\sigma$ over $X_1\cup X_2$ defined by $\sigma.$ For $u,v\in X_1\cup X_2,$ we say that $u\sim_\sigma v$ if and only if one of the following two conditions is satisfied:
\begin{enumerate}
    \item There exists some $k \in \mathbb{Z}_{\ge 0}$ such that $u = \sigma^k(v)$ and $\sigma^j(v)\in X_1\cup X_2$ for all
    $0\le j \le k.$
    \item There exists some $k \in \mathbb{Z}_{\ge 0}$ such that $v = \sigma^k(u)$ and $\sigma^j(u)\in X_1\cup X_2$ for all 
    $0\le j \le k.$
\end{enumerate}
We further distinguish \textit{good} equivalence classes. An equivalence class $E$ is called \textit{good} if there exists some element $x\in E$ and number $k \in \mathbb{N}$ such that $\sigma^k(x)\not \in X_1\cup X_2.$ It trivially holds that the number of good classes is $|X_1\cup X_2| - |S_\sigma| = x_1 + x_2 - \ell(\sigma).$ For each good class $E,$ define by $m(E)$ the unique element $u\in E$ for which $\sigma^{-1}(u)\not \in X_1\cup X_2$ and by $M(E)$ the unique element $v\in E$ for which $\sigma(v)\not \in X_1\cup X_2.$ Both $m(E)$ and $M(E)$ are well-defined since $E$ is a good equivalence class and $\sigma$ is a bijection.\\

\noindent
Now, consider the set $X_0(\sigma)$ of good equivalence classes defined by $\sigma.$ We will construct from $\sigma$ a bijection $\tau$ of $Y(\sigma) = X_0(\sigma)\cup\bigcup_{i = 3}^nX_i$ such that
$X_0(\sigma)\cap\tau (X_0(\sigma)) =\emptyset$ and
$ X_i \cap \tau (X_i)=\emptyset$ for $i \ge 3.$ For $u \in Y(\sigma),$ the image $\tau(u)$ is constructed as follows:
\begin{itemize}
    \item If $u \not \in X_1\cup X_2$ and $\sigma(u)\not \in X_1\cup X_2,$ then $\tau(u): = \sigma(u).$
    \item If $u \not \in X_1\cup X_2$ and $\sigma(u) \in X_1\cup X_2,$ then $\tau(u): = E(u),$ where $E(u)$ is the unique good equivalence $E$ class for which $\sigma(u)  =m(E).$
    \item If $u \in X_0(\sigma),$ then $\tau(u):=\sigma(M(u)).$ 
\end{itemize}

\noindent
Thus far, for every bijection $\sigma$ we have constructed a corresponding bijection over $Y(\sigma).$ We can similarly show how to do the reverse. Namely, given a set $X_0$ of size $x_1+x_2 - \ell$ and bijection $\tau$ of 
$Y = X_0 \cup \bigcup_{i = 3}^nX_i,$ for which $\tau(X_i)\cap X_i = \emptyset$ for all $i,$ we can construct exactly $f(x_1, x_2, \ell)$ bijections $\sigma$ over $X$ for which $X_i \cap \sigma(X_i) = \emptyset$ holds for all $i.$ To do this, we simply need to construct $\sigma$ over $X_1\cup X_2$ such that $|(X_1\cup X_2)\cap\sigma(X_1\cup X_2)| = \ell$ and then revert  the construction with good equivalence classes. By \cref{prop:finterpretation}, this can be done in $f(x_1, x_2, \ell)$ ways. Summing over $\ell$ gives the result.
\end{proof}

\subsection{Schur-Concavity}
The analysis for a fixed number of elements is nearly trivial once we have \cref{thm:recurrence}. To present it, we need the following simple statement.
\begin{proposition}
\label{prop:helpineq1}
If $x_1 \ge x_2\ge 1$ and $\ell \in \mathbb{N},$ then $f(x_1,x_2,\ell)\ge f(x_1+1,x_2-1,\ell).$
\end{proposition}
\begin{proof} To prove the statement, we use the following inequality. Whenever $x_1\ge x_2\ge 1,$ it holds that 
$$\binom{x_1}{\ell} \binom{x_2}{\ell}\ge\binom{x_1+1}{\ell}\binom{x_2-1}{\ell}.$$
This result follows from a simple calculation. In particular, the above inequality is equivalent to
\begin{equation*}
    \begin{split}
       & x_1\cdots(x_1-\ell+1)x_2\cdots(x_2-\ell+1)\ge(x_1+1)\cdots(x_1-l+2)(x_2-1)\cdots(x_2-\ell)\\
      \Longleftrightarrow\; & (x_1-\ell+1)x_2\ge(x_1+1)(x_2-\ell)\\
      \Longleftrightarrow \; & \ell(x_1+1)\ge \ell x_2,
    \end{split}
\end{equation*}
which is trivial.
Applying this inequality twice, we obtain
\begin{equation*}
    \begin{split}
       &\binom{x_1}{\ell_1}\binom{x_2}{\ell_1}(\ell_1)! \binom{x_1}{\ell-\ell_1}\binom{x_2}{\ell-\ell_1}(\ell-\ell_1)!\\
       \ge\; &\binom{x_1+1}{\ell_1}\binom{x_2-1}{\ell_1}(\ell_1)! \binom{x_1+1}{\ell-\ell_1}\binom{x_2-1}{\ell-\ell_1}(\ell-\ell_1)!
    \end{split}
\end{equation*}
Summing over $l$ finishes the proof.
\end{proof}

\noindent
An immediate corollary of \cref{prop:helpineq1} is the following statement.

\begin{corollary} 
\label{cor:oneswitch}
If $x_1\ge x_2\ge 1,$ then $$\derangement(x_1,x_2,\ldots,x_n)\ge \derangement(x_1+1,x_2-1,x_3,\ldots x_n).$$
\end{corollary}
\begin{proof}
Suppose that $x_1\ge x_2\ge 1.$
Using \cref{prop:helpineq1} and \cref{thm:recurrence}, we deduce
\begin{equation*}
    \begin{split}
        & \derangement(x_1,x_2,\ldots,x_n) = 
\sum_{\ell=0}^{x_1+x_2}f(x_1,x_2,\ell) \derangement(x_1+x_2-\ell,x_3,\ldots,x_n)\\
\ge \; & \sum_{\ell=0}^{x_1+x_2}f(x_1+1,x_2-1,\ell) \derangement((x_1+1)+(x_2-1)-\ell,x_3,\ldots,x_n) = 
\derangement(x_1+1,x_2-1,\ldots,x_n) .
    \end{split}
\end{equation*}
\end{proof}

\noindent
By repeatedly applying \cref{cor:oneswitch} as in \cref{obs:schurconvexityfactorials}, we deduce the main result of the current section for $\derangement().$

\begin{theorem}
\label{thm:fixedsumder}
Suppose that $(x_1, x_2, \ldots, x_n)\preceq (y_1, y_2, \ldots, y_n)$ in the sense of \cref{def:majorize}. Then, 
$$
\derangement(x_1, x_2, \ldots, x_n)\ge 
\derangement(y_1, y_2, \ldots, y_n).
$$
\end{theorem}

\noindent
Combining \cref{thm:fixedsumder} and \cref{obs:schurconvexityfactorials}, we also deduce the (weaker) corresponding statement for 
$\anagram().$
\begin{theorem}
\label{thm:fixedsumanar}
Suppose that $(x_1, x_2, \ldots, x_n)\preceq (y_1, y_2, \ldots, y_n)$ in the sense of \cref{def:majorize}. Then, 
$$
\anagram(x_1, x_2, \ldots, x_n)\ge 
\anagram(y_1, y_2, \ldots, y_n).
$$
\end{theorem}

\noindent
Dividing both sides of the inequality in \cref{thm:fixedsumder} by $(x_1 + x_2 + \cdots + x_n)!,$ we also obtain the following statement.

\begin{proposition}
Suppose that $(x_1, x_2, \ldots, x_n)\preceq (y_1, y_2, \ldots, y_n)$ in the sense of \cref{def:majorize}. Then, 
$$
\frac{\anagram(x_1, x_2, \ldots, x_n)}{
\binom{x_1 + x_2 + \cdots + x_n}{x_1, x_2, \ldots, x_n}}\ge 
\frac{\anagram(y_1, y_2, \ldots, y_n)}{
\binom{y_1 + y_2 + \cdots + y_n}{y_1, y_2, \ldots, y_n}
}.
$$
\end{proposition}

\noindent
The last two statement can be rephrased as follows. Words ``lower'' in the poset $\poset_{x_1 + x_2 + \cdots + x_m}$ simultaneously have more anagrams without fixed letters and have a larger fraction of anagrams that have no fixed letters.
\subsection{Beyond Words of the Same Size}
\label{sec:beyondsamesize}
\noindent
We again need statements of the form of \cref{prop:helpineq1}. We defer their proofs to \cref{app:sec4.3} due to their simplicity.

\begin{proposition}
\label{prop:helpineq2}
Suppose that $x_1<x_2.$ Then,
$$
f(x_1, x_2, \ell)<f(x_1+1, x_2, \ell+1).
$$
\end{proposition}

\begin{proposition}
\label{prop:helpineq3}
Suppose that $2x_1<x_2.$ Then,
$$
(x_1 + 1)f(x_1, x_2, \ell)<f(x_1+1, x_2, \ell+1).
$$
\end{proposition}

\noindent
These results are enough to establish the following two theorems.

\begin{theorem}
\label{thm:varsumder}
Suppose that $x_i < \max(x_1, x_2, \ldots, x_n).$ Then,
$$
\derangement(x_1, x_2, \ldots, x_{i-1}, x_i, x_{i+1}, \ldots, x_n)\le
\derangement(x_1, x_2, \ldots, x_{i-1}, x_i+1, x_{i+1}, \ldots, x_n).
$$
\end{theorem}
\begin{proof} Without loss of generality, let $i = 1$ and $2\in \arg \max(x_1, x_2, \ldots, x_n)$ for notational simplicity. Then, using \cref{prop:helpineq2} and 
\cref{thm:recurrence}
, we have 
\begin{equation*}
    \begin{split}
       &  \derangement(x_1+1,x_2,\ldots,x_n)=\sum_{\ell=0}^{x_1+x_2+1}f(x_1+1,x_2,\ell) \derangement(x_1+1+x_2-\ell,\ldots,x_n)\\
       \ge\; & 
       \sum_{\ell=0}^{x_1+x_2}f(x_1+1,x_2,\ell+1) \derangement(x_1+1+x_2-(\ell+1),\ldots,x_n)\\
       \ge \; &
       \sum_{\ell=0}^{x_1+x_2}f(x_1,x_2,\ell) \derangement(x_1+x_2-\ell,\ldots,x_n)=\derangement(x_1,x_2,\ldots,x_n).
    \end{split}
\end{equation*}
\end{proof}

\noindent
As already discussed, the condition $x_i < \max(x_1, x_2, \ldots, x_n)$ is tight. Namely, in the case $x_i = \max(x_1, x_2, \ldots, x_n)$ such an inequality does not always hold. For example $\derangement(t,t)>0 = \derangement(t+1, t)$ for $t \ge 1.$
Unfortunately, we obtain a weaker result - which is likely not tight - for the function $\anagram().$
\begin{theorem}
\label{thm:varsuman}
Suppose that $x_i < \frac{1}{2}\max(x_1, x_2, \ldots, x_n).$ Then,
$$
\anagram(x_1, x_2, \ldots, x_{i-1}, x_i, x_{i+1}, \ldots, x_n)\le
\anagram(x_1, x_2, \ldots, x_{i-1}, x_i+1, x_{i+1}, \ldots, x_n).
$$
\end{theorem}
\begin{proof}
Again, without loss of generality, let $i = 1$ and $2\in \arg \max(x_1, x_2, \ldots, x_n).$ Then, using \cref{prop:helpineq3} and 
\cref{thm:recurrence}
, we have 
\begin{equation*}
    \begin{split}
       &  \derangement(x_1+1,x_2,\ldots,x_n)=\sum_{\ell=0}^{x_1+x_2+1}f(x_1+1,x_2,\ell) \derangement(x_1+1+x_2-\ell,\ldots,x_n)\\
       \ge\; & 
       \sum_{\ell=0}^{x_1+x_2}f(x_1+1,x_2,\ell+1) \derangement(x_1+1+x_2-(\ell+1),\ldots,x_n)\\
       \ge \; &
       \sum_{\ell=0}^{x_1+x_2}(x_1 + 1)f(x_1,x_2,\ell) \derangement(x_1+x_2-\ell,\ldots,x_n)=(x_1 + 1)\derangement(x_1,x_2,\ldots,x_n).
    \end{split}
\end{equation*}
Dividing by $(x_1 + 1)!x_2!\cdots x_n!$ on both sides, we obtain the result.
\end{proof}

\noindent
We leave as an open question the task of improving the constant $\frac{1}{2}.$ Our conjecture is that it can be improved all the way up to 1.

\section{The Anagraph}
\label{se:anagraph}
\subsection{A Universal Reduction of the Alphabet Size}
\label{sec:reduction}
Our main tool in studying global properties of anagraphs will be a technique of simplifying the graph. We achieve this via a reduction of the alphabet size. The reduction works by merging two letters. Specifically, for two distinct letters $i,j,$ and an $n$-tuple $(x_1, x_2, \ldots, x_n),$ we create the $(n-1)$-tuple \linebreak $red_{i,j}(x_1, x_2, \ldots, x_n) = (y_1, y_2, \ldots, y_{n-1}).$ The numbers in $(y_1, y_2, \ldots, y_{n-1})$ are (in some arbitrary order) $x_i + x_j$ and $x_s$ for $s\not \in \{i,j\}.$ 
Intuitively, this reduction makes the letters $i$ and $j$ indistinguishable and everything else remains the same. The following property demonstrates the power of this reduction.

\begin{proposition}
\label{prop:reduction}
For any $n$-tuple $(x_1, x_2, \ldots, x_n)\in \mathbb{Z}_{\ge 0}^{n}$ and distinct letters $i,j\in [n],$ if\linebreak $\anagraphs(red_{i,j}(x_1, x_2, \ldots, x_n))$ satisfies one of the following five properties - 1) Connectivity, 2) Hamiltonicity, 3) Hamilton-Connectivity, 4) Pancyclicity, and 5) Edge-Pancyclicity - then $\anagraphs(x_1, x_2, \ldots, x_n)$ satisfies the same property.
\end{proposition}
\begin{proof}
We will only prove here statement 1) for connectivity  - which we will use in \cref{thm:generalconnectivity} - and  statement 2) for hamiltonicity. The rest of the proofs are very similar and, for that reason, we defer them to \cref{app:sec5.1}.\\

\noindent
Without loss of generality, let $i = n - 1, j = n.$ 
Thus, we will simply take $red_{n-1,n}(x_1, x_2, \ldots, x_n) = (x_1, x_2, \ldots, x_{n-2}, x_{n-1}+x_n).$
We construct the following function $g,$ which takes as input two words $\omega$ and $\lambda,$ where $\omega$ is over $[n-1]$ and has exactly $x_i$ letters $i$ for all $i \le n-2$ and $x_{n-1} + x_n$ letters $n-1,$ and $\lambda$ is a word over $\{n-1, n\}$ and has exactly $x_{n-1}$ letters $n-1$ and $x_{n}$ letters $n.$ Then, $g(\omega, \lambda)$ is the word $\omega$ in which the letters $n-1$ in $\omega$ are replaced with the word $\lambda$ such that their order is preserved. For example, for $n = 4,\omega = 1313323, \lambda = 
3443,
$ we have $g(1313323, 3443) = 1314423.$
Similarly, we define a ``first-argument pseudo-inverse'' of $g,$ denoted by $f.$ It takes as input a word over $[n]$ and replaces each letter $n$ with $n-1.$ For example, when $n = 4,$ $f(1314423) = 1313323.$ More generally, $f(g(\omega, \lambda)) = \omega$ for all choices of $\omega$ and $\lambda.$ Define also a ``second-argument pseudo inverse'' $s$ such that 
$s(g(\omega, \lambda)) = \lambda.$
\\

\noindent
Now, we are ready to present the proof.\\

\noindent
1) \textbf{Connectivity.} Suppose that $\anagraphs(x_1, x_2, \ldots, x_{n-2}, x_{n-1}+x_n)$ is connected. We want to show that for any two words $\chi_1, \chi_2\in V(\anagraphs(x_1, x_2, \ldots, x_n)),$ there is a path in $\anagraphs(x_1, x_2, \ldots, x_n)$ connecting $\chi_1$ and $\chi_2.$ Let $\omega_1 = f(\chi_1),\omega_2 = f(\chi_2),\lambda_1 = s(\chi_1),$ and $\lambda_2 = s(\chi_2).$ Since\linebreak $\anagraphs(x_1, x_2, \ldots, x_{n-2}, x_{n-1}+x_n)$ is connected, there exists a path $\omega_1 = \xi_1, \xi_2, \ldots, \xi_k = \omega_2$ between $\omega_1$ and $\omega_2$ in $\anagraphs(x_1, x_2, \ldots, x_{n-2}, x_{n-1}+x_n).$ One can easily check that
$$
\chi_1 = g(\xi_1, \lambda_1), g(\xi_2, \lambda_1), \ldots, 
g(\xi_{k-1}, \lambda_1), 
g(\xi_k, \lambda_2)  =\chi_2
$$
is a path between $\chi_1$ and $\chi_2$ in $\anagraphs(x_1, x_2, \ldots, x_n).$\\

\noindent
2) \textbf{Hamiltonicity.} Suppose that $\anagraphs(x_1, x_2, \ldots, x_{n-2}, x_{n-1}+x_n)$ is hamiltonian and $\omega_1, \omega_2, \ldots, \omega_k$ is a Hamilton cycle in the graph. Let $\lambda_1, \lambda_2,\ldots, \lambda_\ell$ be all words (in some arbitrary, say alphabetical, order) over the alphabet $\{n-1,n\}$ containing exactly $x_{n-1}$ letters $n-1$ and $x_n$ letters $n.$ Then, clearly 
\begin{equation*}
    \begin{split}
        & g(\omega_1, \lambda_1), g(\omega_2, \lambda_1),
        \ldots, g(\omega_k, \lambda_1),\\
         & g(\omega_1, \lambda_2), g(\omega_2, \lambda_2),
        \ldots, g(\omega_k, \lambda_2),\\
        & \vdots\\
         & g(\omega_1, \lambda_\ell), g(\omega_2, \lambda_\ell),
        \ldots, g(\omega_k, \lambda_\ell)\\
    \end{split}
\end{equation*}
is a hamiltonian cycle in 
$\anagraphs(x_1, x_2, \ldots, x_{n-2}, x_{n-1},x_n).$
\end{proof}

\begin{remark}
\normalfont
\label{rmk:diameter}
The proof for connectivity also demonstrates the following fact. For any $i,j,$ it is the case that 
$$
diam(\anagraphs(x_1,x_2, \ldots,x_n))\le 
diam(\anagraphs(red_{i,j}(x_1,x_2, \ldots,x_n))).
$$
This property, in conjunction with the proof of \cref{lem:3lettercase} can be used to study the diameter of anagraphs.
\end{remark}

\begin{remark}
\normalfont
When we revert the reduction, the five properties are not necessarily preserved. Consider the case $x_1 = x_2 = x_3 = x_4 = 1.$ Then, $\anagraphs(x_1, x_2, x_3, x_4)$ is edge-pancyclic \cite{edgepancyclic}. However, \linebreak $\anagraphs(red_{3,4}(x_1, x_2, x_3, x_4)) = \anagraphs(1,1,2)$ is not even connected (see \cref{thm:generalconnectivity}).
\end{remark}
\subsection{Connectivity of the Anagraph}
\cref{prop:reduction} allows us to fully determine when an anagraph is connected. Specifically, we have the following claim.

\begin{theorem}
\label{thm:generalconnectivity}
Let $(x_1, x_2, \ldots, x_n)\in \mathbb{N}^n.$
The anagraph $\anagraphs(x_1, x_2, \ldots, x_n)$ is connected if and only if one of the following conditions is satisfied:
\begin{enumerate}
    \item $n = 1.$
    \item $n = 2$ and $x_1 = x_2 = 1.$
    \item $n \ge 3, $ the inequality $x_i <\frac{1}{2}\sum_{j = 1}^n x_j$ holds for all indices $i\in [n],$ and $(x_1, x_2, \ldots, x_n)\neq (1,1,1).$
\end{enumerate}
\end{theorem}

\noindent
Before proving this statement in full generality, we will prove two special cases that will be useful in the general proof.

\begin{lemma}
\label{lem:hallmarriageconnectivity}
For $n \ge 4$ and any $k \in \mathbb{N},$ the anagraph 
$\anagraphs(\underbrace{k,k,\ldots, k}_n)$ is connected.
\end{lemma}
\begin{proof}
Let $\chi = \canonicalword(\underbrace{k,k,\ldots, k}_n).$ We will show that for any word $\omega$ in $V(\anagraphs(\underbrace{k,k,\ldots, k}_n)),$ there exists a path from $\chi$ to $\omega$ in the anagraph. This will clearly be enough as anagraphs are vertex-transitive.\\

\noindent
Write $\omega$ as $\chi_1(\omega)\chi_2(\omega)\cdots\chi_n(\omega).$ We will first show the following fact. We can partition the set of positions $[kn]$ into $k$ disjoint sets $S_1, S_2, \ldots, S_k$ which satisfy the following three properties
\begin{enumerate}
    \item $|S_i| = n$ for all $i.$
    \item $\omega|_{S_i}$ contains $n$ different letters for all $i.$
    \item Each $S_i$ contains exactly one position from each interval $[\ell k + 1, (\ell+1)k].$
\end{enumerate}
The proof of this fact follows from Hall's marriage theorem. Consider a bipartite multigraph $\mathcal{G}$ with parts $L$ (stands for letters) and $S$ (stands for subwords). The vertices of $L$ are $[n]$ and the vertices of $S$ are $\{\chi_1(\omega),\chi_2(\omega),\cdots,\chi_n(\omega)\}.$ We connect a subword $\chi_i(\omega)$ to the letter $j$ with multiplicity $t$ if $\chi_i(\omega)$ contains exactly $t$ letters $j.$ Clearly, the resulting graph is $k$-regular. It is by now a folklore fact that the edges of $\mathcal{G}$ can be decomposed into $k$ perfect matchings - $M_1, M_2, \ldots, M_k$ \cite{konig}.\\
Now, arbitrarily match each edge $e = (j, \chi_i(\omega))$ of $\mathcal{G}$ to exactly one position $p(e)$ in $\chi_i(\omega)$ such that $\omega|_{\{p(e)\}} = j$ and each position is matched to exactly one edge. We are ready to form the sets $S_1, S_2, \ldots, S_k.$ We have
$$
S_i = \{p(e)\; : \; e\in M_i\}\text{ for all }i. 
$$

\noindent
We go back to the original problem. Note that for all $i,$ both $\chi|_{S_i}$ and $\omega|_{S_i}$ are permutations of $\lambda =  123\ldots n.$\\

\noindent
However, the derangement graph on $n \ge 4$ vertices is vertex-pancyclic \cite{edgepancyclic}. In particular, the derangement graph is connected and every vertex appears in a cycle of length 3 and a cycle of length 4. As $3$ and $4$ are coprime, this implies that there exists some sufficiently large natural number $N(n)$ such that every two vertices of the derangement graph are connected by a path of length exactly $N(n).$\\

\noindent
Going back to the original problem, we can choose $k$ paths in the derangement graph given by\linebreak $\sigma_{i,1}, \sigma_{i,2}, \ldots, \sigma_{i,N(n)}$ for $1\le i \le k,$ each of length $N(n),$ with the following property.
For each $i,$ it is the case that 
$\sigma_{i,1} = \chi|_{S_i}$ and 
$\sigma_{i,N(n)} = \omega|_{S_i}.$ Now, for all $j \in \{1,2,\ldots, N(n)\},$ define by $\xi_j$ the word with $kn$ letters over alphabet $[n],$ which satisfies that 
$\xi_j|_{S_i} = \sigma_{j,i}$ for all $i \in [k].$ Then, clearly, 
$$
\chi = \xi_1, \xi_2, \ldots, \xi_{N(n)}  =\omega 
$$
is a path in $\anagraphs(\underbrace{k,k,\ldots, k}_n)$ between $\chi$ and $\omega,$ which completes the proof.
\end{proof}

\begin{lemma}
\label{lem:3lettercase}
For $n =3$ and any triplet $(x_1, x_2,x_3)$ of positive integers other than $(1,1,1)$ that satisfies
$x_i < \frac{1}{2}(x_1 + x_2 + x_3)$ for all $i \in \{1,2,3\},$ the anagraph $\anagraphs(x_1,x_2,x_3)$ is connected.
\end{lemma}
\begin{proof}
Consider any word $\omega$ over $\{1,2,3\}$ that has exactly $x_i$ letters $i.$ We will show that for any two positions $u$ and $v$ in $\omega,$ there is a path from $\omega$ to $\omega',$ where $\omega'$ is the same as $\omega,$ except that the letters $\omega|_{\{u\}}$ and $\omega|_{\{v\}}$ are swapped.\\

\noindent
Since the graph is vertex transitive, we only need to show this for the canonical word
$\chi = \canonicalword(x_1,x_2,x_3).$ Without loss of generality, we will show that there exists a path from 
$\chi$ to the word $\lambda$ given by 
$$
\lambda  = 
2\underbrace{11\ldots1}_{x_1-1}
1\underbrace{22\ldots2}_{x_2-1}
\underbrace{33\ldots3}_{x_3}.
$$
We distinguish two cases.\\

\noindent
\textbf{Case 1.} If $x_3 = 1.$ Then, necessarily $x_1 = x_2 = x>1.$ It simply follows that the following path satisfies the desired property.
\begin{equation*}
    \begin{split}
      \chi = \; &  
\underbrace{11\ldots1}_{x}
\underbrace{22\ldots2}_{x}
3\\
&3\underbrace{2\ldots2}_{x-1}
\underbrace{11\ldots1}_{x} 2\\
&2\underbrace{1\ldots1}_{x-1} 
\underbrace{2 \ldots 2}_{x-1} 3 1\\
&1 \underbrace{2 \ldots 2}_{x-1} 
3 \underbrace{1 \ldots 1}_{x-1} 2\\
\lambda = \; &2 \underbrace{1 \ldots 1}_{x-1} 
1 \underbrace{2 \ldots 2}_{x-1} 3.\\
    \end{split}
\end{equation*}

\noindent
\textbf{Case 2.} If $x_3 > 1.$ Let $ u = \min(x_1, x_3-1)>0,$ and $v = x_3 - u>0.$ Clearly, $u + v = x_3 < x_1+x_2$ and $v \le \max(1,x_3-x_1)\le x_2.$
The following path satisfies the desired property. 
\begin{equation*}
    \begin{split}
\chi = \;   
&\underbrace{1 1 \ldots 1}_{x_1} 
\underbrace{2 2 \ldots \;  2}_{x_2} 
\underbrace{3 3 \ldots \;  3}_{x_3} \\
&\underbrace{3 \ldots 3}_{u} 
\underbrace{2 \ldots 2}_{x_1 - u} 
\underbrace{3 \ldots 3}_{v} 
\underbrace{1 \ldots 1}_{x_2 - v} 
\underbrace{2 \ldots 2}_{x_2 -x_1+ u} 
\underbrace{1 \ldots 1}_{x_1 -x_2 + v} \\
\lambda = \; 
&2 \underbrace{1 \ldots \; 1}_{x_1-1} 
1 \underbrace{2 \ldots \;  2}_{x_2-1} 
\underbrace{3 3 \ldots \;  3}_{x_3}.\\
    \end{split}
\end{equation*}

\end{proof}

\noindent
Now, we are ready to handle the general case.

\begin{proof}[Proof of \cref{thm:generalconnectivity}]
We consider four cases based on $n.$\\

\noindent
\textbf{Case 1.} $n = 1.$ Then, the anagraph has a single vertex, no matter what $x_1$ is.\\

\noindent
\textbf{Case 2.} $n = 2.$ If $x_1 = x_2 = 1,$ then the anagraph has only two vertices $(1,2)$ and $(2,1)$ and is, therefore, connected. If $x_1 \neq x_2,$ then the anagraph is clearly disconnected as it has more than vertex, but all of its vertices are isolated (see \cref{obs:positivity}). If $x_1 = x_2 = x>1,$ consider the canonical word $\chi = 
\underbrace{1\ldots 1}_{x}
\underbrace{2\ldots 2}_{x}.
$ It is clearly only connected to $\omega = \underbrace{2\ldots 2}_{x}
\underbrace{1\ldots 1}_{x},$ but $\omega$ also has no other neighbours than $\chi.$ In particular, $\chi$ is in a different connected component from $\lambda = 2\underbrace{1\ldots 1}_{x-1}
1\underbrace{2\ldots 2}_{x-1}.$\\

\noindent
\textbf{Case 3.} First, we know that if $x_i < \frac{1}{2}(x_1 + x_2 + x_3)$ holds for all $i \in \{1,2,3\}$ and $(x_1, x_2, x_3)\neq (1,1,1),$ the anagraph is connected by \cref{lem:3lettercase}. Now, we need to show that no other anagraph in this case is connected. First, $\anagraphs(1,1,1)$ is not connected as it has two connected components - $\{123,231,312\}$ and $\{132,321,213\}.$ Now, we need to show that if $x_i \ge \frac{1}{2}(x_1 + x_2 + x_3)$ holds for some $i,$ the anagraph is also disconnected. Without loss of generality, let this be the case for $i = 1.$ If 
$x_1 > \frac{1}{2}(x_1 + x_2 + x_3),$ then the anagraph has at least one vertex, but all of its vertices are isolated by \cref{obs:positivity}, so the graph is disconnected. If $x_1 = \frac{1}{2}(x_1 + x_2 + x_3),$ consider the canonical word $\chi.$ It is clearly only connected to words whose last $x_1$ letters equal $1.$ However, any word whose last $x_1$ letters equal 1 is only connected to words whose first $x_1$ letters equal 1. In particular, this means that for any word $\omega$ in the connected component of $\chi,$ either its first $x_1$ letters or its last $x_1$ letters are all $1$'s. Thus, the anagraph is disconnected in this case as well.\\

\noindent
\textbf{Case 4.} $n \ge 4.$ First, if there exists some $i$ such that $x_i \ge \frac{1}{2}\sum_{j = 1}^nx_j,$ the anagraph is disconnected. This follows in the same way as in the case $n  = 3.$\\
If, on the other hand, $x_i < \frac{1}{2}\sum_{j = 1}^nx_j$ holds for all $i,$ we need to show that the anagraph is connected. We do so by repeatedly reducing the alphabet size as follows.\\
Let $(y^{(n)}_1, y^{(n)}_2, \ldots, y^{(n)}_n) = (x_1, x_2, \ldots, x_n).$ While $k >4,$ find $u,v$ such that $y^{(k)}_u$ and $y^{(k)}_v$ are the smallest two (tiebreaks handled arbitrarily) numbers in $(y^{(k)}_1, y^{(k)}_2, \ldots, y^{(k)}_k).$ Then, define \linebreak
$(y^{(k-1)}_1, \ldots, y^{(k-1)}_{k-1}) = red_{u,v}(y^{(k)}_1, y^{(k)}_2, \ldots, y^{(k)}_k).$\\
Note that in doing this operation, the sum $\sum_{j = 1}^{k}y^{(k)}_j$ remains unchanged. This shows that the condition $y^{(k)}_i< \frac{1}{2}\sum_{j = 1}^k y^{(k)}_j$ holds for all $i$ and $k\ge 4.$ Indeed, the only number that appears in $(y^{(k)}_1, \ldots, y^{(k)}_{k}),$ but not in $(y^{(k+1)}_1, y^{(k+1)}_2, \ldots, y^{(k+1)}_{k+1}),$ is $y^{(k+1)}_u + y^{(k+1)}_u.$ As $k \ge 4,$ the choice of $u,v,$ guarantees that  
$y^{(k+1)}_u + y^{(k+1)}_v\le \frac{2}{5}\sum_{j = 1}^{k}y^{(k)}_j< \frac{1}{2}\sum_{j = 1}^{k}y^{(k)}_j.$\\

\noindent
After performing this operation, we are left with four numbers $y^{(4)}_1, y^{(4)}_2, y^{(4)}_3, y^{(4)}_4.$ Write them in decreasing order as $y_1\ge y_2\ge y_3\ge y_4.$ We know that $y_1 \le \frac{1}{2}(y_1 + y_2 + y_3 + y_4).$ We now consider two cases for these numbers.\\

\noindent
\textbf{Case 4.1.} If $y_1 = y_2 = y_3 = y_4,$ then $\anagraphs(y_1,y_2,y_3,y_4)$ is connected by \cref{lem:hallmarriageconnectivity}. From \cref{prop:reduction}, so is $\anagraphs(x_1, x_2, \ldots, x_n).$\\

\noindent
\textbf{Case 4.2.} If the four numbers are not equal, this means that $y_3 + y_4 <y_1 + y_2.$ Therefore, \linebreak
$\anagraphs(y_1,y_2,y_3+y_4)$ is connected by \cref{lem:3lettercase}.  From \cref{prop:reduction}, so is $\anagraphs(x_1, x_2, \ldots, x_n).$
\end{proof}

\begin{remark}
\normalfont
We end with a remark about a further global property of the anagraphs, beyond the ones listed in \cref{def:hamiltonicity}. Combining together our result for the parity of the degrees in an anagraph - \cref{cor:parity} - and our result for the connectivity of anagraphs - \cref{thm:generalconnectivity} - we derive a necessary and sufficient condition for an anagraph to be Eulerian.
\end{remark}

\section{Further Directions}
As the non-asymptotic study of anagrams without fixed letters is a rather new field, we end with a multitude of directions for further work. Of course, we are most interested in fully resolving
\cref{prob:arithmetic}, \cref{prob:ordinal}, and \cref{prob:hamiltonicityproblem}. Here, we make a few remarks about these problems.\\

\noindent
\textbf{Arithmetic Questions:} We are especially interested in improving the running-time of \cref{cor:algorithmprime}, so that larger primes can also be handled efficiently. As discussed in \cref{rmk:productofprimes}, this direction of study can be very useful in computing $\anagram().$ Similarly, one could think of extending the result to prime-powers.
\\

\noindent
\textbf{Ordinal Questions:} Since our result about Schur-Concavity makes very substantial progress in the case of $\sum_i x_i = \sum_j y_j,$ we believe that 
the setting of equal sums is more tractable. Beyond that, we are interested in improving the constant $\frac{1}{2}$ in \cref{thm:varsuman}. Namely, we make the following conjecture.

\begin{conjecture}
Suppose that $x_i <\max(x_1, x_2, \ldots, x_n).$ Then,
$$
\anagram(x_1, x_2, \ldots, x_{i-1}, x_i, x_{i+1}, \ldots, x_n)\le
\anagram(x_1, x_2, \ldots, x_{i-1}, x_i+1, x_{i+1}, \ldots, x_n).
$$
\end{conjecture}
  
\noindent
\textbf{Questions on Anagraphs:} In light of the strong connection between hamiltonicity and vertex-transitivity \cite{transitive}, we make the following conjecture. 

\begin{conjecture}
All, but finitely many, connected anagraphs are hamiltonian.
\end{conjecture}

\noindent
We exclude finitely many anagraphs to avoid $\anagraphs(1,1)$ (which is a single edge) and the four known connected vertex-transitive graphs which are not hamiltonian \cite{edgepancyclic}. Note that if this conjecture is true, one potentially only needs to prove it for anagraphs over alphabet of size four. This follows from \cref{prop:reduction} and the method used in the proof of \cref{thm:generalconnectivity}.

\section*{Acknowledgements}
I want to thank Aleksander Makelov and Lyuben Lichev for the many helpful suggestions and comments they made during the development of this work when I first started it in 2018. I am also grateful to Noah Kravitz who helped me improve the paper by pointing out a hole in the proof of \cref{thm:primemodulo} in an earlier version. Finally, I am thankful to Ilia Bozhinov for his assistance with a computer implementation for computing $\anagram(),$ which was especially helpful while I was working on the arithmetic aspects of the function. This research did not receive any specific grant from funding agencies in the public, commercial, or not-for-profit sectors.

\printbibliography
\appendix
\section{Omitted Proofs in Section 4.3}
\label{app:sec4.3}
\begin{proof}[Proof of \cref{prop:helpineq2}]
We simply need to show that for any $\ell\le 2x_1$ and $0\le \ell_1\le \ell,$ it is the case that 
\begin{equation*}
    \begin{split}
        & \binom{x_1+1}{\ell_1+1}\binom{x_2}{\ell_1+1}(\ell_1+1)! \binom{x_1+1}{\ell-\ell_1}\binom{x_2}{\ell-\ell_1}(\ell-\ell_1)!\\
        \ge \; &
        \binom{x_1}{\ell_1}\binom{x_2}{\ell_1}(\ell_1)! \binom{x_1}{\ell-\ell_1}\binom{x_2}{\ell-\ell_1}(\ell-\ell_1)!\\
    \end{split}
\end{equation*}
Clearly, we only need to consider the case when $\ell_1\le x_1$ as otherwise both sides equal $0.$ In this case, we open the parenthesis as follows
\begin{equation*}
    \begin{split}
         &(x_1+1)\cdots(x_1-\ell_1+1)x_2 \cdots(x_2-\ell_1)(x_1+1)\cdots(x_1-(\ell-\ell_1)+2)x_2\cdots(x_2-(\ell-\ell_1)+1)\\
        \ge \; &
        (\ell_1+1)x_1\cdots(x_1-\ell_1+1)x_2\cdots(x_2-\ell_1+1)x_1\cdots(x_1-(\ell-\ell_1)+1)x_2\cdots(x_2-(\ell-\ell_1)+1)\\
        \Longleftrightarrow\; &
        (x_1+1)(x_1+1)(x_2-\ell_1)\ge(\ell_1+1)(x_1-(\ell-\ell_1)+1).
    \end{split}
\end{equation*}
The last inequality holds since $\ell_1\le \ell, \ell_1\le x_1, \ell_1+1\le x_2.$ Indeed, this implies that 
$$(\ell_1+1)(x_1-(\ell-\ell_1)+1)\le(x_1+1)(x_1+1)\le(x_1+1)(x_1+1)(x_2-\ell_1).$$
\end{proof}

\begin{proof}[Proof of \cref{prop:helpineq3}] The proof is almost the same. We want to show that under the given conditions, 
\begin{equation*}
    \begin{split}
        & \binom{x_1+1}{\ell_1+1}\binom{x_2}{\ell_1+1}(\ell_1+1)! \binom{x_1+1}{\ell-\ell_1}\binom{x_2}{\ell-\ell_1}(\ell-\ell_1)!\\
        \ge \; &
        (x_1 + 1)\binom{x_1}{\ell_1}\binom{x_2}{\ell_1}(\ell_1)! \binom{x_1}{\ell-\ell_1}\binom{x_2}{\ell-\ell_1}(\ell-\ell_1)!\\
    \end{split}
\end{equation*}
When we open the brackets, this reduces to 
$$
(x_1+1)(x_2-\ell_1)\ge(\ell_1+1)(x_1-(\ell-\ell_1)+1).
$$
The last inequality holds as $x_1 <\frac{1}{2}x_2,\ell_1\le \ell,$ and $\ell_1\le x_1.$ Indeed, this implies that
$$
(x_1+1)(x_2-\ell_1)\ge (x_1 + 1)(2x_1 + 1-x_1)  =
(x_1 + 1)(x_1 + 1)\ge 
(\ell_1+1)(x_1-(\ell-\ell_1)+1).
$$
\end{proof}

\section{Omitted Proofs in Section 5.1}
\label{app:sec5.1}
\begin{proof}[Omitted Proofs in \cref{prop:reduction}]
1) \textbf{Hamilton-Connectivity.} Suppose that\linebreak $\anagraphs(x_1, x_2, \ldots, x_{n-2}, x_{n-1}+x_n)$ is hamilton-connected. In particular, this means that it is hamiltonian since we can choose an edge $(\omega,\xi)\in E(\anagraphs(x_1, x_2, \ldots, x_{n-2}, x_{n-1}+x_n))$ and, together with the hamiltonian path connecting $\omega$ and $\xi,$ this edge will form a hamilton cycle.\\ 
We want to show that for any two words $\chi_1, \chi_2\in V(\anagraphs(x_1, x_2, \ldots, x_n)),$ there is a hamilton path in $\anagraphs(x_1, x_2, \ldots, x_n)$ connecting $\chi_1$ and $\chi_2.$ Let $\omega_1 = f(\chi_1),\omega_2 = f(\chi_2).$ Since\linebreak $\anagraphs(x_1, x_2, \ldots, x_{n-2}, x_{n-1}+x_n)$ is hamilton connected, there exists a hamilton path $\omega_1 = \xi_1, \xi_2, \ldots,$ $\xi_k = \omega_2$ between $\omega_1$ and $\omega_2$ in $\anagraphs(x_1, x_2, \ldots, x_{n-2}, x_{n-1}+x_n).$ 
Similarly, there exists a hamilton cycle 
$\omega_1 = \zeta_1, \zeta_2, \ldots, \zeta_k$ in 
$\anagraphs(x_1, x_2, \ldots, x_{n-2}, x_{n-1}+x_n).$ 
Let $\lambda_1, \lambda_2,\ldots, \lambda_\ell$ be all words (the order will be determined later) over the alphabet $\{n-1,n\}$ containing exactly $x_{n-1}$ letters $n-1$ and $x_n$ letters $n.$ We now distinguish two cases.\\
\textbf{Case 1)}
If
$s(\chi_1) \neq s(\chi_2).$ Without loss of generality, we can index $\lambda_1, \lambda_2, \ldots, \lambda_\ell$ such that $s(\chi_1) = \lambda_1, s(\chi_2) = \lambda_n.$ Then, the following is a hamilton path connecting the two vertices:
\begin{equation*}
    \begin{split}
        & g(\zeta_1, \lambda_1), g(\zeta_2, \lambda_1),
        \ldots, g(\zeta_k, \lambda_1),\\
         & g(\zeta_1, \lambda_2), g(\zeta_2, \lambda_2),
        \ldots, g(\zeta_k, \lambda_2),\\
        & \vdots\\
         & g(\zeta_1, \lambda_{\ell-1}), g(\zeta_2, \lambda_{\ell-1}),
        \ldots, g(\zeta_k, \lambda_{\ell-1}),\\
        & g(\xi_1, \lambda_{\ell}), g(\xi_2, \lambda_{\ell}),
        \ldots, g(\xi_k, \lambda_{\ell})
    \end{split}
\end{equation*}
\textbf{Case 2)} If $s(\chi_1) = s(\chi_2).$ 
Without loss of generality, we can index $\lambda_1, \lambda_2, \ldots, \lambda_\ell$ such that $s(\chi_1) = \lambda_1.$ Then, 
the following is a hamilton path connecting the two vertices:
\begin{equation*}
    \begin{split}
        & g(\zeta_1, \lambda_1), g(\zeta_2, \lambda_2),
        g(\zeta_3, \lambda_2),
        \ldots, g(\zeta_k, \lambda_2),\\
         & g(\zeta_1, \lambda_2), g(\zeta_2, \lambda_3),
         g(\zeta_3, \lambda_3),
        \ldots, g(\zeta_k, \lambda_3),\\
        & \vdots\\
         & g(\zeta_1, \lambda_{\ell-1}), g(\zeta_2, \lambda_{\ell}),
         g(\zeta_3, \lambda_{\ell}),
        \ldots, g(\zeta_k, \lambda_{\ell}),\\
        & g(\xi_1, \lambda_{\ell}), g(\xi_2, \lambda_{1}),
        g(\xi_3, \lambda_{1}),
        \ldots, g(\xi_k, \lambda_{1}).
    \end{split}
\end{equation*}

\noindent
4) \textbf{Pancyclicity.} Suppose that $\anagraphs(x_1, x_2, \ldots, x_{n-2}, x_{n-1}+x_n)$ is pancyclic. In particular, this means that it is hamiltonian. Furthermore, as the graph is vertex-transitive, it is also \textit{vertex-pancyclic}, which means that for any $\omega$ in its vertex set and number $3\le k\le 
|V(\anagraphs(x_1, x_2, \ldots, x_{n-2}, x_{n-1}+x_n))|,
$ 
there exists a simple cycle of length $k$ containing the vertex $\omega.$\\

\noindent
We want to show that 
$\anagraphs(x_1, x_2, \ldots, x_{n-2}, x_{n-1},x_n)$
is also pancyclic.
Take any $k$ such that 
\begin{equation*}
    \begin{split}
         3\le k \le \; &\\
\le \; & |V(\anagraphs(x_1, x_2, \ldots,x_n))|\\
=\;  &\binom{x_1 + x_2 + \cdots + x_n}{x_1, x_2, \ldots, x_{n-2}, x_{n-1},x_n} \\
=\; & \binom{x_1 + x_2 + \cdots + x_n}{x_1, x_2, \ldots, x_{n-2}, x_{n-1}+x_n}\binom{x_{n-1}+x_n}{x_n}\\
=\;  &  
|V(\anagraphs(x_1, x_2, \ldots, x_{n-2}, x_{n-1}+x_n))|
\binom{x_{n-1}+x_n}{x_n}. 
    \end{split}
\end{equation*}
We need to show that there exists a cycle in $\anagraphs(x_1, x_2, \ldots,x_n)$ of size $k.$ We assume that\linebreak $k< |V(\anagraphs(x_1, x_2, \ldots,x_n))|$ as we have already shown that the reduction preserves hamiltonicity.\\
Denote $s = |V(\anagraphs(x_1, x_2, \ldots, x_{n-2}, x_{n-1}+x_n))|.$
Suppose that $k = m s + r,$ where $r$ is the residue upon division by $s$ and $m< \binom{x_{n-1}+x_n}{x_n}$ is the quotient.  Since $s\ge 5 $ (one can trivially check that there does not exist a pancyclic anagraph on less than 5 vertices with \cref{thm:generalconnectivity}), we can clearly write $k$ as a sum of $t\in \{m,m+1\}$ numbers $s_1, s_2, \ldots, s_t,$ such that 
$3 \le s_i \le s$ for all $i\in [t].$ Since $t \le \binom{x_{n-1}+x_n}{x_n},$ we can also choose $t$ different words $\lambda_1, \lambda_2, \ldots, \lambda_t,$ each composed of exactly $x_{n-1}$ letters $n-1$ and $x_n$ letters $n.$ Furthermore, for each $i \in [t],$ there exists a cycle $\xi_{i,1}, \xi_{i,2}, \ldots, \xi_{i,s_i} $ such that $\xi_{i,s_i} = \xi_{i+1, 1}$ when $i <t$ also holds. This follows from the fact that the reduced anagraph is vertex-pancyclic. Thus,
\begin{equation*}
    \begin{split}
         & g(\xi_{1,1}, \lambda_1), 
         g(\xi_{1,2}, \lambda_1),\ldots,
        g(\xi_{1,s_1}, \lambda_1),\\
        & g(\xi_{2,1}, \lambda_2), 
         g(\xi_{2,2}, \lambda_2),\ldots,
        g(\xi_{2,s_2}, \lambda_2),\\
        \vdots\\
        & g(\xi_{t,1}, \lambda_t), 
         g(\xi_{t,2}, \lambda_t),\ldots,
        g(\xi_{t,s_t}, \lambda_t)\\
    \end{split}
\end{equation*}
is a simple cycle of length $k$ in 
$\anagraphs(x_1, x_2, \ldots,x_n).$\\

\noindent
5) \textbf{Edge-Pancyclicity.} 
Edge pancyclicity follows in absolutely the same way as pancyclicity, except that we choose the first cycle
$\xi_{1,1}, \xi_{1,2}, \ldots, \xi_{1,s_1}$ so that 
$(\xi_{1,1},\xi_{1,2})$ corresponds to the desired edge
$(\chi_1, \chi_2)$
of $\anagraphs(x_1, x_2, \ldots,x_n)$ to be included. 
We distinguish two cases $s(\chi_1) = s(\chi_2)$ and  $s(\chi_1) \neq s(\chi_2)$ and handle them in absolutely the same way as for hamilton-connectivity. 
\end{proof}
\end{document}